\documentclass[11pt,a4paper,leqno]{amsart}
\usepackage{amsmath} 
\usepackage{amssymb} 
\usepackage{amscd}
\usepackage{amsxtra}
\newtheorem{theorem}{Theorem}[section]
\newtheorem{proposition}[theorem]{Proposition}
\newtheorem{corollary}[theorem]{Corollary}
\newtheorem{lemma}[theorem]{Lemma}

\newtheorem{remark}[theorem]{Remark}

\newtheorem{conjecture}[theorem]{Conjecture}

\newcommand{\C}{{\mathbb C}}
\newcommand{\R}{{\mathbb R}}
\newcommand{\Z}{{\mathbb Z}}
\newcommand{\T}{{\mathbb T}}
\newcommand{\F}{{\mathbb F}}
\newcommand{\Proj}{{\mathbb P}}

\newcommand{\bs}{{\boldsymbol s}}
\newcommand{\bm}{{\boldsymbol m}}
\newcommand{\bd}{{\boldsymbol d}}
\newcommand{\bc}{{\boldsymbol c}}
\newcommand{\be}{{\boldsymbol e}}
\newcommand{\bL}{\mathbb L}
\newcommand{\vx}{{\vec{x}}}

\newcommand{\mH}{{\mathcal{H}}}
\newcommand{\mL}{{\mathcal{L}}}
\newcommand{\mO}{{\mathcal{O}}}
\newcommand{\mI}{{\mathcal{I}}}
\newcommand{\mIas}{\mI^{\rm asym}}
\newcommand{\mP}{{\mathcal{P}}}
\newcommand{\mPT}{\mathcal{P}^\T}
\newcommand{\mR}{{\mathcal{R}}}
\newcommand{\mT}{{\mathcal{T}_{\mM}}}
\newcommand{\mM}{{\mathcal{M}}}
\newcommand{\mD}{{\mathcal{D}}}
\newcommand{\mC}{{\mathcal{C}}}

\newcommand{\frm}{{\mathfrak{m}}}
\newcommand{\td}{w}          
\newcommand{\tdD}{W}         
\newcommand{\mc}{\mathsf{w}} 
\newcommand{\mcT}{\mathsf{T}} 
\newcommand{\mct}{\mathsf{t}}
\newcommand{\mcx}{\mathsf{x}}
\newcommand{\cp}{\mathsf{p}} 
\newcommand{\cu}{\mathsf{u}} 
\newcommand{\ncu}{{\tt u}}   

\newcommand{\tI}{{\widetilde{I}}}
\newcommand{\tnabla}{\widetilde{\nabla}^\hbar}
\newcommand{\tL}{\widetilde{L}}
\newcommand{\tT}{\widetilde{T}}
\newcommand{\FH}{{FH_{S^1}^*}}
\newcommand{\tFH}{{\widetilde{FH}_{S^1}^*}}
\newcommand{\FHT}{{FH_{\T\times S^1}^*}}
\newcommand{\DeltaT}{\Delta^\T}

\newcommand{\hq}{\hat{q}}
\newcommand{\hatt}{\hat{t}}
\newcommand{\hx}{\hat{x}}
\newcommand{\hatotimes}{\widehat{\otimes}}

\newcommand{\hess}{\triangle}
\newcommand{\tEu}{\widetilde{E}} 
\newcommand{\hOmega}{\widehat{\Omega}}
\newcommand{\homega}{\widehat{\omega}}
\newcommand{\ovd}{{\overline{d}}}

\newcommand{\Spec}{\operatorname{Spec}}
\newcommand{\id}{\operatorname{id}}
\newcommand{\End}{\operatorname{End}}

\newcommand{\Mat}{\operatorname{Mat}}
\newcommand{\Res}{\operatorname{Res}}
\newcommand{\ad}{\operatorname{ad}}
\newcommand{\ch}{\operatorname{ch}}
\newcommand{\Ch}{\operatorname{Ch}}
\newcommand{\Euler}{\operatorname{Euler}} 
\newcommand{\emb}{\operatorname{\mathfrak{emb}}}
\newcommand{\Crit}{\operatorname{Crit}}
\newcommand{\crit}{\operatorname{crit}}
\newcommand{\Vol}{\operatorname{Vol}}
\newcommand{\Conv}{\operatorname{Conv}}
\newcommand{\Graph}{\operatorname{Graph}}
\newcommand{\Hess}{\operatorname{Hess}}

\def\ov#1{\overline{#1}}
\def\corr#1{\left\langle{#1}\right\rangle}
\def\corrL#1{\left\langle{#1}\right\rangle^{\mathcal{L}}}
\def\parti#1{{q_{#1}\frac{\partial}{\partial q_{#1}}}}
\def\parfrac#1#2{{\frac{\partial #1}{\partial #2}}}
\def\pair#1#2{\langle #1,#2\rangle}
\def\pairL#1#2{\langle #1,#2\rangle^{\mathcal{L}}}
\def\pairT#1#2{\left\langle #1,#2 \right\rangle_{\T}}
\def\pairr#1#2{\langle\!\langle #1,#2\rangle\!\rangle}

\begin{document}

\title{Convergence of Quantum Cohomology by Quantum Lefschetz}
\author{Hiroshi Iritani}
\address{Department of Mathematics, Graduate School of Science, 
Kyoto University, Oiwake-cho, Kitashirakawa, Sakyo-ku, Kyoto, 606-8502, Japan.}
\email{iritani@math.kyoto-u.ac.jp}
\begin{abstract}
Quantum Lefschetz theorem by Coates and Givental \cite{coates-givental} 
gives a relationship between the genus 0 Gromov-Witten theory of $X$ 
and the twisted theory by a line bundle $\mL$ on $X$. 
We prove the convergence of the twisted theory under the assumption 
that the genus 0 theory for original $X$ converges. 
As a byproduct, we prove the semisimplicity and the Virasoro conjecture for 
the Gromov-Witten theories of (not necessarily Fano) 
projective toric manifolds. 
\end{abstract}
\maketitle 

\section{Introduction}

Quantum cohomology is a 
deformation of the ring structure of the ordinary cohomology.  
The structure constants of quantum cohomology are  
formal power series whose coefficients consist of  
Gromov-Witten invariants. 
We do not know a priori whether or not 
the structure constants are convergent. 
In this paper, we discuss the compatibility 
of quantum Lefschetz principle and the convergence of quantum cohomology. 

There are several cases where the convergence is trivial. 
If $c_1(X)>0$, the small quantum cohomology of $X$ is defined over 
the polynomial ring by the degree constraints. 
If $c_1(X)<0$, even the big quantum cohomology is defined 
over the polynomial ring for the same reason. 
Hence, the problem is the intermediate case, i.e. 
when there exist two curves $C_1, C_2$ in $X$ such that 
$\pair{c_1(X)}{[C_1]}\ge 0$ and $\pair{c_1(X)}{[C_2]}\le 0$. 
The main theorem in this paper is the following. 
\begin{theorem}
\label{thm:mainthm}
Let $X$ be a smooth projective variety and $\mL$ be a nef line bundle 
on $X$. If the big quantum cohomology $QH^*(X)$ of $X$ 
has convergent structure constants, then the twisted quantum cohomology 
$QH^*_{S^1}(X,\mL)$ by $\mL$ also has convergent structure constants. 
(This holds true when $\mL$ is replaced by a sum of nef line bundles. )
\end{theorem}

Coates-Givental's quantum Lefschetz theorem \cite{coates-givental} 
gives the twisted quantum cohomology $QH^*_{S^1}(X,\mL)$ 
in terms of $QH^*(X)$. 
Here, $QH^*_{S^1}(X,\mL)$ is a cohomology theory 
closely related to the quantum cohomology of an intersection $Y\subset X$
with respect to the line bundle $\mL$.  
More precisely, $QH^*_{S^1}(X,\mL)$ gives us 
the information on the structure constants of $QH^*(Y)$ 
with respect to the cohomology classes coming from the ambient space $X$. 
Therefore, if the convergence of $QH^*(X)$ is known, 
we can know the convergence of $QH^*(Y)$ partially. 
The main tool in the proof is a ring of formal power series 
with certain estimates for coefficients.

In the second half of the paper, 
we give a description of mirror symmetry 
for a not necessarily nef toric variety. 
In \cite{iritani-genmir}, 
the author calculated the quantum cohomology $D$-module of a toric variety $X$. 
The method there was to embed $X$ into another Fano toric variety $X'$ 
as a complete intersection and 
to use quantum Lefschetz theorem \cite{coates-givental} 
together with a mirror theorem for a Fano toric variety \cite{givental-mirrorthm-toric}. 
We will recast the consequences of \cite{iritani-genmir} 
in terms of the following oscillatory integral, 
which was introduced as a mirror of a toric variety
in \cite{givental-ICM,givental-mirrorthm-toric}: 
\[
\mI_\Gamma(q_1,\dots,q_r,\hbar)=\int_{\Gamma_q\subset Y_q}
e^{\sum_{i=1}^{r+N}\mc_i/\hbar} \omega_q, \quad 
Y_q=\left\{(\mc_i)_{i=1}^{r+N}\in (\C^*)^{r+N}\; ;\; 
\textstyle\prod_{i=1}^{r+N} \mc_i^{m_{ia}}=q_a\right \}.   
\]
Here, $\omega_q$ is a holomorphic volume form on $Y_q$ and 
$\Gamma_q$ is a non-compact cycle. 
These oscillatory integrals define 
the following mirror $D$-module $M_{\rm mir}$ 
(denoted by $FH_0$ in the main text): 
\[
M_{\rm mir}:=\C\langle q_1,\dots,q_r,
\hbar\partial_1,\dots,\hbar\partial_r,\hbar\rangle 
/I_{\rm poly}, \quad \partial_a=q_a\partial/\partial q_a.
\]
Here, $I_{\rm poly}$ is a left ideal consisting of polynomial
differential operators annihilating $\mI_{\Gamma}(q,\hbar)$. 
Mirror symmetry for a Fano toric variety $X$ 
states that the mirror $D$-module $M_{\rm mir}$ is isomorphic to 
the big quantum cohomology $D$-module $QDM^*(X)$ restricted to $H^2(X)$. 
Here, $QDM^*(X)$ is a $D$-module over the total cohomology ring $H^*(X)$
which is defined by $QH^*(X)$.   
For a general toric variety $X$, we obtain the following description:  

\begin{theorem}
[see Theorem \ref{thm:generalizedmirrortransformation} for details]
Let $\widehat{M}_{\rm mir}$ be a $D$-module on the formal germ
$(\C^r,0)$ obtained as 
the completion of $M_{\rm mir}$ with respect to its natural $q$-adic topology 
(denoted by $\FH$ in the main text). 
There exists a formal embedding $\emb\colon (\C^r,0)\rightarrow (H^*(X),0)$ 
such that we have an isomorphism of $D$-modules: 
\[
\Phi_{\emb}\colon \emb^*(QDM^*(X)) \cong \widehat{M}_{\rm mir}. 
\]
Here, $QDM^*(X)$ is the big quantum $D$-module of the toric variety $X$.   
If $X$ is Fano, the image of $\emb$ coincides with 
the linear subspace $H^2(X)\subset H^*(X)$. 
\end{theorem}

Because of the completion in the above description, 
it is not clear if $QH^*(X)$ is convergent. 
Our main theorem \ref{thm:mainthm} is not directly applicable to $X$  
because $X$ is a complete intersection in $X'$ 
with respect to a sum of not necessarily nef line bundles.
Using techniques similar to the proof of Theorem \ref{thm:mainthm}, however, 
we show the following: 
\begin{theorem}[Theorem \ref{thm:analyticityofgenmirtrans},
Corollary \ref{cor:semisimplicity}]
The big quantum cohomology of a smooth projective toric variety
is convergent and generically semisimple. 
The embedding $\emb$ in the above theorem is complex analytic.  
\end{theorem} 

Note that the isomorphism $\Phi_{\emb}$ is not convergent unless $X$ is nef 
(Proposition \ref{prop:nefcharacterization}). 
The asymptotic expansion of the oscillatory integral 
$\mI_{\Gamma}(q,\hbar)$ in $\hbar$ is shown to 
give a formal solution to $QDM^*(X)$ for {\it special} choices of cycles $\Gamma$
(Corollary \ref{cor:asymptoticsolution}). 
We also prove the $R$-conjecture 
for equivariant quantum cohomology of toric varieties 
(Theorem \ref{thm:R-conjecture}).  
Here, the $R$-conjecture implies the Virasoro constraints 
by Givental's theory \cite{givental-quadratic}.

Our result on the semisimplicity is also a successful test 
for Bayer and Manin's modified Dubrovin's conjecture \cite{bayer-manin}.  
The modified Dubrovin's conjecture claims that 
$(p,p)$-part of quantum cohomology of a projective variety $X$ 
is generically semisimple if 
its bounded derived category $D^b_{\rm coh}(X)$ of coherent sheaves 
admits a full exceptional collection.  
In fact, Kawamata recently showed that toric varieties have  
full exceptional collections \cite{kawamata}. 

For the application of the main theorem \ref{thm:mainthm}, 
we need to know the convergence of big quantum cohomology of ambient spaces. 
We prove that if $H^*(X)$ is generated by $H^2(X)$ 
and if the small quantum cohomology of $X$ 
has convergent structure constants, so does the big quantum cohomology 
(Corollary \ref{cor:convergenceofbig}). 
In particular, the big quantum cohomology of a Fano variety 
with $H^2$-generated cohomology  
always has convergent structure constants. 
In a subsequent paper \cite{iritani-localization}, 
we will also prove that 
the big quantum cohomology (and higher genus potential 
$\mathcal{F}_g$ also) 
is convergent for a projective manifold 
which admits Hamiltonian torus action with  
only isolated fixed points and isolated one dimensional orbits.

We should remark that in this paper, 
we only consider the even part of (quantum) cohomology. 
For example, $H^*(X)$ always means $H^{\rm even}(X)$.

The paper is organized as follows.
In section 2 and 3, we review the quantum $D$-modules 
and the quantum Lefschetz theorem by Coates and Givental. 
In section 4, we prove the main theorem. 
In section 5, we discuss the mirror symmetry for a non-nef toric variety. 
In section 6, we prove the $R$-conjecture
for any toric variety. 

\noindent
{\bf Acknowledgments.}
Thanks are due to Professor Martin Guest, Professor Hiraku Nakajima 
and Kazushi Ueda for valuable discussions. 
Part of this paper was written 
while the author stayed at Mathematical Sciences Research Institute.  
He thanks MSRI for excellent working conditions.  
He is also grateful to anonymous referee for thier helpful comments.  
This research is supported by Grant-in-Aid for JSPS Fellows 
and Scientific Research 15-5482.

\section{Quantum $D$-modules}
\label{sect:quantumD-modules}

In this section, we introduce quantum $D$-modules 
twisted by the equivariant Euler class following 
\cite{coates-givental, pandharipande}.
Let $X$ be a smooth projective variety and $\mL$ be a line bundle over $X$. 
Let $\ov{M}_{0,n}(X,\bd)$ be the moduli space of genus 0, degree $\bd$ 
stable maps to $X$ with $n$ marked points, where $\bd\in H_2(X,\Z)$. 
We have the following diagram: 
\begin{equation*}
\begin{CD} 
\ov{M}_{0,n+1}(X,\bd) @>{e_{n+1}}>>  X \\
@V{\pi_{n+1}}VV  \\
\ov{M}_{0,n}(X,\bd)
\end{CD}
\end{equation*}
where $e_i$ is the evaluation map and $\pi_i$ is the forgetful map. 
We introduce the fiber-wise $S^1$ action on $\mL$ by scalar multiplication. 
Let $\lambda$ be a generator of the $S^1$ equivariant cohomology of a point. 
Define the twisted correlator 
\[\corrL{\alpha_1,\dots,\alpha_n}_{S^1,\bd}
 =\int_{[\ov{M}_{0,n}(X,\bd)]^{\rm virt}}
 \prod_{i=1}^n e_i^*(\alpha_i)\Euler_{S^1}
(R^\bullet{\pi_{n+1}}_*e_{n+1}^*\mL),  
\]
where $\alpha_1,\dots,\alpha_n\in H^*(X,\C)$
and  $[\ov{M}_{0,n}(X,\bd)]^{\rm virt}$ is the virtual fundamental class. 
The right hand side is in $\C[\lambda,\lambda^{-1}]$. 
Let $\{p_0,\dots,p_s\}$ be a basis of $H^*(X,\C)$. 
We assume that $p_0$ is a unit and that 
$p_1,\dots ,p_r$ form a nef integral basis of $H^2(X,\Z)$ 
($r\le s$). 
Let $t_0,\dots,t_s$ be linear coordinates dual to the basis 
$p_0,\dots,p_s$.  
We write $q_a:= \exp(t_a)$ and 
$q^\bd:=q_1^{\pair{p_1}{\bd}} q_2^{\pair{p_2}{\bd}}
\cdots q_r^{\pair{p_r}{\bd}}=\exp(\sum_{a=1}^r \pair{p_a}{\bd}t_a)$ 
for $\bd\in H_2(X,\Z)$.  
Define a twisted pairing $\pairL{\cdot}{\cdot}_{S^1}$ by 
\[
\pairL{\alpha}{\beta}_{S^1}=\int_{X} \alpha\cup\beta\cup \Euler_{S^1}(\mL).
\]
The twisted quantum product $*_\mL$ is defined by the formula: 
\begin{align*}
\pairL{\alpha*_\mL\beta}{\gamma}_{S^1}
 & =\sum_{\bd\in \Lambda}\sum_{n\ge 0}\frac{1}{n!}
    \corrL{\alpha,\beta,\gamma,
       (\textstyle\sum_{j=0}^s t_j p_j)^{\otimes n}}_{S^1,\bd} \\ 
 &= \sum_{\bd\in \Lambda}\sum_{n\ge 0}\frac{1}{n!}
    \corrL{\alpha,\beta,\gamma,
       (\textstyle\sum_{j=r+1}^s t_j p_j)^{\otimes n}}_{S^1,\bd} q^\bd,
\end{align*}
where $\alpha,\beta,\gamma \in H^*(X,\C)$ and 
$\Lambda\subset H_2(X,\Z)$ is a semigroup generated by effective curves. 
The product $*_\mL$ is extended linearly over 
$\C[\lambda,\lambda^{-1}][\![t_0,q_1,\dots,q_r,t_{r+1},\dots,t_s]\!]$ 
and 
\[
QH^*_{S^1}(X,\mL):=(H^*(X)\otimes
\C[\lambda,\lambda^{-1}][\![t_0,q_1,\dots,q_r,t_{r+1},\dots,t_s ]\!],*_\mL) 
\]
becomes an associative and commutative ring. 
We define $\deg t_j=2-\deg p_j$ for $j=0$ or $j>r$, 
$\deg q^\bd=2\pair{c_1(X)-c_1(\mL)}{\bd}$ and $\deg \lambda=2$. 
Then $QH^*_{S^1}(X,\mL)$ becomes a graded ring. 
For simplicity, we also use the following notation: 
\[ 
x=(x_0,x_1,\dots,x_r,x_{r+1},\dots,x_s)=(t_0,q_1,\dots,q_r,t_{r+1},\dots, t_s),\quad 
\tau=\sum_{j=0}^s t_j p_j.
\]
When $c_1(\mL)$ is nef, 
the product $*_\mL$ can be defined over $\C[\lambda][\![x]\!]$ 
(see \cite{pandharipande}). 
Therefore in this case, 
we can also consider the non-equivariant ($\lambda=0$) version
$QH^*(X,\mL)=(H^*(X)\otimes \C[\![x]\!], *_\mL)$.
  
The usual (non-twisted) quantum cohomology 
$QH^*(X)=(H^*(X)\otimes\C[\![x]\!],*_X)$ is defined by removing all the Euler classes
in the definition.  
Let $Y$ be a smooth intersection in $X$ with respect to 
a transverse section $s\in \Gamma(X,\mL)$. 
By the main theorem of \cite{kim-kresch-pantev}, if $c_1(\mL)$ is nef, 
we have at $\lambda=0$, 
\[
\corrL{\alpha_1,\dots,\alpha_n}_\bd=
\corr{i^*\alpha_1,\dots,i^*\alpha_n}^Y_\bd
\]
for the inclusion $i\colon Y\rightarrow X$  
and the correlator $\corr{\cdots}^Y_\bd$ for $Y$.   
Therefore, from $*_\mL$, we can read the structure constants of $QH^*(Y)$ 
with respect to the classes coming from the ambient space. 

We can endow a $D$-module structure on the quantum cohomology. 
Let $*$ denote $*_\mL$ or $*_X$. 
The dual Givental connection $\nabla^\hbar$ is defined by 
\begin{gather*}
\nabla^\hbar_j = \hbar \parfrac{}{t_j} + p_j* \quad (0\le j\le s),
\end{gather*}
where $\hbar$ is a formal variable of degree two. 
This connection is regular singular along $q_1=\dots=q_r=0$ 
and is known to be flat. 
It defines the non-twisted or twisted 
quantum $D$-modules: 
\[ QDM^*(X)= (H^*(X)[\hbar][\![x]\!],\nabla^\hbar),\quad 
QDM^*_{S^1}(X,\mL)= (H^*(X)[\hbar,\lambda,\lambda^{-1}][\![x]\!],\nabla^\hbar)
\]
When $c_1(L)$ is nef, we can also consider the non-equivariant 
version $QDM^*(X,\mL)$. 
It is easy to see that quantum $D$-module is generated 
by $1$ over the Heisenberg algebra 
$\C[\hbar][\![x]\!][\nabla^\hbar_0,\dots,\nabla^\hbar_s]$.  
There exists a unique fundamental solution $L(\tau,\hbar)$ for the 
flat connection $\nabla^\hbar$ such that 
\begin{gather}
\label{eq:fundamentalsolutionL}
\hbar d \circ L(\tau,\hbar)= L(\tau,\hbar) \circ \nabla^\hbar, \\
\nonumber
L(\tau,\hbar)=e^{\tau/\hbar}T(x,\hbar), \quad 
T\in \End(H^*(X))\otimes \C[\hbar^{-1}][\![x]\!], \quad T|_{q=0}=\id.  
\end{gather}
where $\tau\in H^*(X)$ is considered as an operator 
acting on $H^*(X)$ by the cup product.
This $L(\tau,\hbar)$ satisfies the following unitarity
\cite{givental-elliptic}: 
\begin{equation}
\label{eq:unitarity}
\pair{L(\tau,-\hbar)\alpha}{L(\tau,\hbar)\beta}=\pair{\alpha}{\beta}   
\end{equation}
where $\pair{\cdot}{\cdot}$ denotes the Poincar\'{e} pairing 
$\pair{\cdot}{\cdot}^X$ of $X$ in case of $QDM^*(X)$ 
and the twisted Poincar\'{e} pairing $\pairL{\cdot}{\cdot}_{S^1}$ 
in case of $QDM^*_{S^1}(X,\mL)$. 
This $L(\tau,\hbar)$ also defines the $J$-function of quantum cohomology 
by  
\[
J(\tau,\hbar):=L(\tau,\hbar) 1, \quad 
1 \in H^*(X) \text{ is a unit.}
\]
The $J$-function is a realization of a generator $1$ of the quantum $D$-module 
as a function.  

\section{Coates-Givental's Quantum Lefschetz}
\label{sect:quantumLefschetz}
Coates and Givental's quantum Lefschetz theorem \cite{coates-givental}
describes the relationship between the two quantum cohomologies  
$QH^*(X)$ and $QH^*_{S^1}(X,\mL)$. 
It was described in terms of symplectic transformations of 
Lagrangian cones in the infinite dimensional 
space $H^*(X)\otimes\C[\![\hbar,\hbar^{-1}]\!]$. 
In this paper, we describe it in terms of gauge transformations
by translating the language of cones into that of quantum $D$-modules.

\subsection{Symplectic formalism}
First we will review the infinite dimensional symplectic formalism 
in \cite{givental-quadratic,coates-givental} briefly. 
Consider the following general multiplicative characteristic class 
for a complex vector bundle $E$: 
\[
\bc(E)=\exp\left(\sum_{k=0}^\infty s_k \ch_k(E)\right).  
\] 
Here, $\bs=(s_0,s_1,s_2,\dots)$ are infinite number of arbitrary parameters.  
Let $\mH$ be the following infinite dimensional vector space: 
\[
\mH:=H^*(X)\otimes\C[\hbar,\hbar^{-1}][\![\bs]\!]=\mH_+\oplus \mH_-, 
\]
where 
\[
\mH_+=H^*(X)\otimes\C[\hbar][\![\bs]\!], \quad 
\mH_-=\hbar^{-1}H^*(X)\otimes\C[\hbar^{-1}][\![\bs]\!]. 
\]
For any holomorphic vector bundle $E$ on $X$, 
define a symplectic form $\Omega_\bs$ on $\mH$ by 
\[
\Omega_\bs(f(\hbar),g(\hbar))=\Res_{\hbar=0} d\hbar
\int_X f(-\hbar) g(\hbar) \bc(E),\quad f,g\in \mH.
\]
Then, $\mH_+$ and $\mH_-$ become Lagrangian with respect to this symplectic 
form  $\Omega_\bs$ and give a polarization of $\mH$. 
The fundamental solution $L(\tau,\hbar)$ of the quantum $D$-module $QDM^*(X)$ 
(see (\ref{eq:fundamentalsolutionL})) defines a 
Lagrangian subspace $\bL_\tau$ of $(\mH,\Omega_0)$ for each $\tau$. 
\begin{equation}
\label{eq:Lagrangiansubspace}
\bL_\tau:= L(\tau,-\hbar)(\mH_+) 
\end{equation}
Here, we assume the convergence of $L(\tau,\hbar)$ 
for the sake of simplicity. 
For a rigorous argument, we need to introduce a Novikov ring $\Lambda_{\rm nov}$ 
and replace $\mH$ with the module 
$H^*(X,\Lambda_{\rm nov}[\![\bs]\!])\{\hbar,\hbar^{-1}\}$ 
of convergent power series in $\hbar$ with 
respect to the adic topology (as explained in \cite{coates-givental}). 
These semi-infinite subspaces sweep a germ of Lagrangian cone $\mC_0$ 
in $(\mH,\Omega_0)$. 
\begin{equation}
\label{eq:Lagrangiancone}
\mC_0:=\bigcup_{\tau\in H^*(X)} \hbar\bL_\tau 
\end{equation}
As explained in \cite{coates-givental}, 
the tangent space of $\mC_0$ at any point in $\hbar\bL_\tau$
equals the Lagrangian subspace $\bL_\tau$. 
In other words, $\mC_0$ has a remarkable property that 
{\it it is ruled by $\hbar$ times its tangent spaces}. 
Let $J(\tau,\hbar)=L(\tau,\hbar) 1$ be the $J$-function of $QH^*(X)$.  
Then the vectors $-\hbar J(\tau,-\hbar)$ parametrized by $\tau\in H^*(X)$ 
lie on the cone $\mC_0$. 
The derivatives $\{-\hbar \partial_{t_j}J(\tau,-\hbar)\}_j$ form 
a basis of the tangent space $\bL_\tau$ of $\mC_0$ 
over $\C[\hbar][\![\bs]\!]$. 
Thus, the $J$-function recovers the whole Lagrangian cone $\mC_0$ 
by (\ref{eq:Lagrangiancone}). 

Similarly, the twisted theory by the characteristic class $\bc$
and a vector bundle $E$ defines a Lagrangian cone 
$\mC_\bs \subset (\mH,\Omega_\bs)$. 
Coates and Givental proved that two cones $\mC_0$ and $\mC_\bs$ are 
related by a linear symplectic transformation.  
\begin{theorem}[\cite{coates-givental}, Corollary 4]
\label{thm:quantumLefschetzoriginal}
The linear symplectic transformation 
\[
\bc(E)^{-1/2}\exp \left(\sum_{l,k\ge 0}
s_{2k+l-1}\frac{B_{2k}}{(2k)!} \ch_l(E)\hbar^{2k-1}\right) 
\colon (\mH,\Omega_0) \longrightarrow (\mH,\Omega_\bs)
\]
sends the Lagrangian cone $\mC_0$ to $\mC_\bs$. 
Here, $B_{2k}$ is the Bernoulli number defined by 
$x/(1-e^{-x})=x/2+\sum_{k\ge 0} B_{2k}x^{2k}/(2k)!$ and $s_{-1}=0$. 
\end{theorem}

In this paper, 
we only consider a twist by the equivariant Euler class
and the case where $E$ is a line bundle $\mL$. 
In this case, values of the parameters $s_i$ are set as follows: 
\begin{equation}
\label{eq:substitution}
s_0=\log \lambda,\quad 
s_k=(-1)^{k-1}\frac{(k-1)!}{\lambda^k} \quad k>0. 
\end{equation}
Using this substitution, $\bc(\mL)$ equals the equivariant Euler class 
$c_1(\mL)+\lambda$, where $S^1$ acts on $\mL$ by scalar multiplication on each fiber 
and $\lambda$ is a generator of $H_{S^1}^*({\rm pt})$.  
Coates and Givental introduced the following {\it hypergeometric modification} 
of the $J$-function: 
\[
I_{\mL}(\tau,\hbar):=
\sum_{\bd\in \Lambda} 
\frac{\prod_{k=-\infty}^{\pair{\rho}{\bd}}(\rho+k\hbar+\lambda)} 
     {\prod_{k=-\infty}^{0}(\rho+k\hbar+\lambda)}
J_\bd(t,\hbar)q^\bd, \quad \rho:=c_1(\mL),
\] 
where $J(\tau,\hbar)=\sum_{\bd\in \Lambda} J_\bd(t,\hbar)q^\bd$ 
is the $J$-function of $QH^*(X)$.  
\begin{theorem}[{\cite[Theorem 2]{coates-givental}}]
\label{thm:hypergeometricmodification} 
The vectors $-\hbar I_{\mL}(\tau^*,-\hbar)$ parametrized by $\tau^*\in H^*(X)$ 
lie on the cone $\mC_\lambda$ defined by the twisted quantum cohomology 
$QH_{S^1}^*(X,\mL)$. 
The derivatives $\{-\hbar\partial_{t_j}I_{\mL}(\tau^*,-\hbar)\}_j$ span the 
tangent spaces to $\mC_\lambda$ and 
reconstruct the whole cone $\mC_\lambda$. 
\end{theorem}


\subsection{Symplectic transformation as a gauge transformation}
Following \cite[section 5]{iritani-genmir}, 
we will interpret the linear symplectic transformation in 
Theorem \ref{thm:quantumLefschetzoriginal}  
as a gauge transformation of quantum $D$-modules. 
The quantum $D$-module 
$QDM^*(X)=(H^*(X)\otimes\C[\hbar][\![x]\!],\nabla^\hbar)$ 
is a vector bundle on a formal neighborhood $\Spec\C[\![x_0,\dots,x_s]\!]$ 
with fiber $H^*(X)$ endowed with a flat connection $\nabla^\hbar$ 
with parameter $\hbar$. 
This vector bundle has a canonical trivialization by definition. 
This canonical trivialization together with a canonical origin $\tau=0$ of 
the cohomology fixes a choice of a fundamental solution 
(\ref{eq:fundamentalsolutionL}) and in turn defines 
the Lagrangian subspace (\ref{eq:Lagrangiansubspace}) and 
the Lagrangian cone (\ref{eq:Lagrangiancone}). 
Suppose that one changes the trivialization by a gauge transformation 
$g(x,\hbar)\in \End(H^*(X))\otimes \C[\hbar][\![x]\!]$ such that 
\[
g_0(\hbar):=g|_{q=0}=A_0\exp(a_1\hbar+a_3\hbar^3+\cdots), \quad A_0, a_i\in H^*(X).
\]  
Then the dual Givental connection $\nabla^\hbar$ changes as 
\[
g^*\nabla^\hbar_j =g^{-1}\circ\nabla^\hbar_j \circ g 
= \hbar \parfrac{}{t_j} + g^{-1}(p_j*)g+\hbar g^{-1}\parfrac{g}{t_j}
\]
Suppose also that one takes another point $\tau=c$ as an origin. 
Then the fundamental solution changes as 
\begin{equation}
\label{eq:changeoffundamentalsolution}
L^{g,c}(\tau,\hbar)=e^{-c/\hbar}g_0(\hbar)^{-1}L(\tau,\hbar) g(x,\hbar). 
\end{equation} 
This new fundamental solution $L^{g,c}(\tau,\hbar)$ is uniquely determined 
by a differential equation 
$\hbar d \circ L^{g,c}(\tau,\hbar)=L^{g,c}(\tau,\hbar)\circ (g^*\nabla^\hbar)$ and 
a shifted initial condition:
\[
L^{g,c}(\tau,\hbar)=e^{(\tau-c)/\hbar} T^{g,c}(x,\hbar), \quad T^{g,c}|_{q=0}=\id.  
\]
By the same formulas (\ref{eq:Lagrangiansubspace}), (\ref{eq:Lagrangiancone}) 
as before,  $L^{g,c}(\tau,\hbar)$ defines the following Lagrangian cone: 
\[
e^{c/\hbar}g_0(-\hbar)^{-1}\mC_0= 
A_0^{-1}\exp\left (\frac{c}{\hbar}+a_1\hbar +a_3\hbar^3+\cdots
\right ) \mC_0. 
\] 
Therefore, $e^{c/\hbar}g_0(-\hbar)^{-1}$ is identified  with the symplectic 
transformation in Theorem \ref{thm:quantumLefschetzoriginal}. 
In case of equivariant Euler class, 
Theorem \ref{thm:quantumLefschetzoriginal} is restated as follows: 
\begin{proposition}
\label{prop:symplectictransformationasagaugetransformation} 
The dual Givental connection of the twisted quantum $D$-module 
$QDM^*_{S^1}(X,\mL)$ is obtained from that of $QDM^*(X)$ by a gauge transformation  
$\nabla^\hbar\mapsto g^*\nabla^\hbar=g^{-1} \circ \nabla^\hbar \circ g $, 
$g(x,\hbar)\in \sqrt{\lambda}\End(H^*(X))
\otimes\C[\hbar][\lambda,\lambda^{-1}]\!][\![x]\!]$ 
and a coordinate change $\tau\mapsto \hat\tau=\hat\tau(\tau)$ 
such that 
\begin{gather}
\label{eq:initialcondofgauge}
g|_{q=0}= (\lambda+\rho)^{1/2}\exp \left(\sum_{k\ge 1,l\ge 0}(-1)^l
\frac{(2k+l-2)!}{\lambda^{2k+l-1}} 
\frac{B_{2k}}{(2k)!} \frac{\rho^l}{l!} \hbar^{2k-1}\right). \\
\label{eq:coordinateshift}
\hat\tau(\tau=c)=0, \quad \text{where } 
c=\rho \log \lambda+
\sum_{l=2}^{\dim X} (-1)^{l}\frac{(l-2)!}{\lambda^{l-1}}\frac{\rho^l}{l!}.
\end{gather}
\end{proposition} 

The hypergeometric modification 
in Theorem \ref{thm:hypergeometricmodification} can be considered as 
an intermediate step to find a gauge transformation $g$ and 
a coordinate change $\hat\tau=\hat\tau(\tau)$ in the above proposition. 
If one changes a trivialization of $QDM^*(X)$ 
by a gauge transformation $g$ satisfying (\ref{eq:initialcondofgauge})
and also shifts the origin by a coordinate change satisfying 
(\ref{eq:coordinateshift}), 
the set of column vectors 
$L^{g,c}(\tau,-\hbar)p_j$ of the new fundamental solution 
gives a basis of tangent spaces to the 
cone $\mC_\lambda$ defined by the twisted quantum cohomology. 
Conversely, a matrix formed by the basis 
$\{-\hbar\partial_{t_j}I_{\mL}(\tau^*,-\hbar)\}_j$
of tangent spaces to the cone $\mC_\lambda$ 
(given in Theorem \ref{thm:hypergeometricmodification}) 
satisfies the following initial condition: 
\begin{equation}
\label{eq:fundamentalsolutiontL}
\tL(\tau^*,\hbar):=
\begin{bmatrix}
\vert & &\vert  \\
\hbar\partial_{t_0}I_{\mL}(\tau^*,\hbar) & \cdots & 
\hbar\partial_{t_s}I_{\mL}(\tau^*,\hbar) \\
\vert&  &\vert  
\end{bmatrix}
=e^{\tau^*/\hbar} \tT(x,\hbar), \quad \tT(x,\hbar)|_{q=0}=\id.    
\end{equation}
From this it follows that there exists 
an intermediate gauge transformation $g_1(\tau,\hbar)$
satisfying (\ref{eq:initialcondofgauge}) and a 
coordinate change $\tau^*=\tau^*(\tau)$ 
satisfying (\ref{eq:coordinateshift}) such that 
\[
L^{g_1,c}(\tau,\hbar)= \tL(\tau^*,\hbar), 
\]
where the left hand side is given in (\ref{eq:changeoffundamentalsolution}). 
Let $\tnabla$ be the connection in this new gauge and coordinates: 
\begin{equation}
\label{eq:tnabla}
\tnabla:=\hbar d+ \sum_{j=0}^s \Omega_j dt_j, \quad 
\Omega_j(x,\hbar):=\tL(\tau,\hbar)^{-1}\hbar \partial_{t_j} 
\tL(\tau,\hbar).
\end{equation} 
Here, the connection $\tnabla$ at a point $\tau^*(\tau)$ 
is equal to $g_1^*\nabla^\hbar$ at $\tau$. 
By abuse of notation, we use $\tau$ instead of $\tau^*$ 
as an argument of this new connection $\tnabla$ and solution $\tL$. 

\begin{proposition}
\label{prop:regularity}
{\rm (i)}  
The connection matrix $\Omega_j$ is in the ring 
$\End(H^*(X))\otimes\C[\hbar][\lambda,\lambda^{-1}]\!][\![x]\!]$. 
If $c_1(\mL)$ is nef, this is also in the ring 
$\End(H^*(X))\otimes\C[\hbar,\lambda][\![x]\!]$. 

{\rm (ii)} The pairing $\pairr{\alpha}{\beta}$ defined by 
\[
\pairr{\alpha}{\beta}:=\pairL{\tL(\tau,-\hbar)\alpha}{\tL(\tau,\hbar)\beta}_{S^1}, 
\quad \alpha,\beta\in H^*(X)
\]
takes values in $\C[\hbar][\lambda,\lambda^{-1}]\!][\![x]\!]$. 
If $c_1(\mL)$ is nef, this takes values in 
$\C[\hbar,\lambda][\![x]\!]$. 
\end{proposition}
\begin{proof}
(i) A tangent space to the cone $\mC_\lambda$ 
is a vector space over $\C[\hbar][\lambda,\lambda^{-1}]\!]$. 
The former part follows from this and the argument in \cite{coates-givental}
using a ruling property (\ref{eq:Lagrangiancone}) of the cone $\mC_\lambda$. 
If $c_1(\mL)$ is nef, the hypergeometric modification is of the form
\[
I_{\mL}(\tau,\hbar) = \sum_{\bd\in \Lambda} 
\prod_{k=1}^{\pair{\rho}{\bd}} (\rho+\lambda+k\hbar) 
J_{\bd}(t,\hbar) q^\bd
\]
and does not contain negative powers of $\lambda$. 
The latter part follows from this. 

(ii) This follows from that $\tL(\tau,-\hbar)\alpha$ is a tangent vector of 
$\mC_\lambda$ and that $\mC_\lambda$ is Lagrangian 
with respect to the symplectic form 
$\Res_{\hbar=0}\pairL{f(-\hbar)}{g(\hbar)}_{S^1} d\hbar$. 
\end{proof}

The dual Givental connection of $QDM^*_{S^1}(X,\mL)$ 
can be obtained from $\tnabla$ by a further gauge transformation by 
$g_2\in \End(H^*(X))\otimes \C[\hbar][\lambda,\lambda^{-1}]\!][\![x]\!]$ 
and a coordinate change 
$x\mapsto \hx$. 
The gauge transformation $g_2$ must satisfy  
\begin{equation}
\label{eq:secondgaugetransformation}
g_2|_{q=0}=\id,\quad 
g_2^*\tnabla =\hbar d +\sum_{j=0}^s \hOmega_j dt_j, 
\quad  \hOmega_j \text{ does not depend on } \hbar.  
\end{equation}
The new coordinates 
$\hx=(\hatt_0,\hq_1=e^{\hatt_1},\dots,\hq_r=e^{\hatt_r},
\hatt_{r+1},\dots,\hatt_s)$ 
are of the form  
\begin{gather}
\label{eq:secondcoordinatechange}
\hatt_0=t_0+F_0(x,\lambda), \quad \log \hq_a=\log q_a+ F_a(x,\lambda) 
\quad (1\le a\le r), \\ 
\nonumber
\hatt_j=t_j+F_j(x,\lambda), \quad (r+1\le j\le s)
\end{gather}
for some $F_j(x,\lambda)\in \C[\lambda,\lambda^{-1}]\!][\![x]\!]$ 
satisfying $F_j(x,\lambda)|_{q=0}=0$. 
When written in the new coordinate system $(\hatt_0,\dots,\hatt_s)$, 
the connection matrix $\hOmega_{\hat{j}}$ must satisfy  
\begin{equation}
\label{eq:secondcoordinatechangecond}
\hOmega_{\hat{j}}(1)=p_j, \quad 
\hOmega_{\hat{j}}:=\sum_{i=0}^s \parfrac{t_i}{\hatt_j} \hOmega_i
\end{equation}
since $\hOmega_{\hat{j}}$ is identified with the quantum multiplication by 
$p_j$ in $QH^*_{S^1}(X,\mL)$. 
Note that a gauge transformation and a coordinate change 
satisfying (\ref{eq:secondgaugetransformation}), (\ref{eq:secondcoordinatechange}) 
do not change the cone $\mC_\lambda$.  
As shown in \cite{coates-givental,guest}, 
this gauge transformation $g_2$ can be obtained as the positive part 
of the Birkhoff factorization of the fundamental solution:
\[
g_2=\tL_+, \quad \tL(\tau,\hbar)=\tL_-(\tau,\hbar) \tL_+(x,\hbar),   
\] 
where $\tL_-=\id + O(\hbar^{-1})$ and $\tL_+$ is regular at $\hbar=0$. 
As shown in \cite{iritani-genmir}, Theorem 4.6 and 4.8, 
$g_2$ and $\hx$ can be uniquely determined 
by the conditions (\ref{eq:secondgaugetransformation}),  
(\ref{eq:secondcoordinatechangecond}). 

Summarizing, we can find the gauge transformation $g$ in Proposition 
\ref{prop:symplectictransformationasagaugetransformation} 
as a composite $g_1\circ g_2$ of two gauge transformations. 
First gauge transformation $g_1$ corresponds to
the hypergeometric modification $\tL(\tau,\hbar)$ of the fundamental solution 
and the second one $g_2$ can be obtained by the Birkhoff factorization 
of $\tL(\tau,\hbar)$. 
In the next section, we study the analytic property of 
the connection $\Omega_j$ and the second 
gauge transformation $g_2$ in detail.

\section{Proof of the main theorem}
\label{sect:proofofmainthm}
\subsection{Formal power series with estimates}
\label{subsect:powerserieswithest}
In this subsection, we give an estimate for the connection matrix 
$\Omega_j$ which is obtained after the gauge transformation $g_1$. 
We will show that matrix elements of $\Omega_j$ are 
not necessarily convergent power series,  
but their coefficients satisfy certain estimates. 
We assume that $c_1(\mL)=\rho$ is nef and 
that the ambient quantum cohomology $QH^*(X)$ converges.
We use the notation in Section 3. 

Let $L(\tau,\hbar)$ be the fundamental solution 
(\ref{eq:fundamentalsolutionL}) of $QDM^*(X)$. 
\[
L(\tau,\hbar)=
e^{\tau/\hbar}T(x,\hbar)=
e^{\tau/\hbar}\sum_{\bm,n\ge 0}T_{\bm,n}x^\bm \hbar^{-n}, 
\]
where $\bm=(m_1,\dots,m_s)$ is a multi-index in $\Z^s_{\ge 0}$. 
($T(x,\hbar)$ does not depend on $x_0=t_0$.)

\begin{lemma}
\label{lem:Jest}
There exist positive constants $C_1, C_2>1$ such that 
\[ 
\|T_{\bm,n}\| \le C_1 C_2^{|\bm|+n}\frac{1}{n!},
\]
where $\|\cdot \|$ is the operator norm and $|\bm|=\sum_{j=1}^s m_j$. 
\end{lemma}
\begin{proof}
The function $T(x,\hbar)$
is known to satisfy the following homogeneity:  
\[
2\hbar\parfrac{T}{\hbar}+
\sum_{a=1}^r (\deg q_a)\parfrac{T}{q_a}+
\sum_{j=r+1}^s(\deg t_j)t_j\parfrac{T}{t_j}+[\mu,T]=0, 
\]
where $\mu$ is a constant matrix defined by $\mu(p_j)=(\deg p_j) p_j$ 
and $\deg q_a$ is the degree defined in non-twisted theory.  
By using $\hbar\partial L/\partial t_j=L(p_j*_X)$, 
we can rewrite the above equation as 
\[
2\hbar\parfrac{T}{\hbar}+
\sum_{a=1}^r (\deg q_a)\frac{1}{\hbar}(T(p_a*_X)-p_a T)+
\sum_{j=r+1}^s (\deg t_j) \frac{t_j}{\hbar}T(p_j*_X)+[\mu,T]=0
\]
Set $T_n(x)=\sum_{\bm}T_{\bm,n}x^\bm$. 
By this equation, we have 
\[
(2n-\ad(\mu))T_n(x)=T_{n-1}(x)A(x)-c_1(X)T_{n-1}(x),
\] 
where  
$A(x)=\sum_{a=1}^r (\deg q_a) (p_a*_X)
+\sum_{j=r+1}^s t_j(\deg t_j)(p_j*_X)$. 
Therefore, we have 
\[
T_n=\frac{1}{2n}\left(1-\frac{\ad(\mu)}{2n}\right)^{-1}
(T_{n-1}A-c_1(X)T_{n-1}).
\]
By the assumption that $QH^*(X)$ is convergent, 
there exist a neighborhood $U$ of $x=0$ and a constant $C>0$ such that 
\[
\sup_{x\in U} \|A(x)\|\le C,\quad 
\sup_{n\ge 1}
\left\|\left(1-\frac{\ad(\mu)}{2n}\right)^{-1}\right\|\le C,\quad
\|c_1(X)\|\le C.
\] 
Since $T_0(x)=\id$, we can see that  
$\sup_{x\in U}\|T_n(x)\|\le C^{2n}/n!$. 
The lemma follows from this estimate.
\end{proof}

For a multi-index $\bm=(m_1,\dots,m_s)$, we write 
$\pair{\rho}{\bm}=\sum_{a=1}^r \xi_a m_a$, 
where $\rho=c_1(\mL)=\sum_{a=1}^r \xi_a p_a$. 
Note that $\pair{\rho}{\bm}\ge 0$ for $\bm\in \Z_{\ge 0}^s$ 
because $\rho$ is nef. 
Let $\tL(\tau,\hbar)$ be the fundamental solution given by hypergeometric modification 
(\ref{eq:fundamentalsolutiontL}). 
\[
\tL(\tau,\hbar)=e^{\tau/\hbar}\tT(x,\hbar,\lambda)
=e^{\tau/\hbar}\sum_{\bm\ge 0,n\in \Z}
\tT_{\bm,n}(\lambda) x^\bm\hbar^{-n}
\]
\begin{lemma}
\label{lem:estofmodification}
There exists a positive constant $C_3(\lambda)$ 
depending continuously on $\lambda$ such that 
\[
\|\tT_{\bm,n}(\lambda)\|\le C_1C_3(\lambda)^{|\bm|+|n|}
\begin{cases}
1/n! & n\ge 0, \\
|n|! & n\le 0.
\end{cases}
\]
Here we set $0!=1$. 
Moreover, $\tT_{\bm,n}=0$ for $-n>\pair{\rho}{\bm}$. 
\end{lemma}
\begin{proof}
Because $\rho$ is nef, we have 
\begin{equation}
\label{eq:hypmodifforT}
\sum_{n}\tT_{\bm,n}(\lambda)\hbar^{-n}=
\prod_{k=1}^{\pair{\rho}{\bm}}(\rho+\lambda+k\hbar)
\sum_{n\ge 0}T_{\bm,n}\hbar^{-n}. 
\end{equation}
Therefore, by Lemma \ref{lem:Jest}, we have 
\begin{align*}
\|\tT_{\bm,n}(\lambda)\| &=
\left\|\sum_{0\le l\le \pair{\rho}{\bm}}
\sum_{1\le k_1<k_2<\cdots<k_l\le\pair{\rho}{\bm}} 
(\rho+\lambda)^{\pair{\rho}{\bm}-l}k_1k_2\cdots k_l T_{\bm,n+l}\right\|\\
&\le
\sum_{0\le l\le \pair{\rho}{\bm}}
\|(\rho+\lambda)\|^{\pair{\rho}{\bm}-l}\binom{\pair{\rho}{\bm}}{l}
(\pair{\rho}{\bm}-l+1)\cdots \pair{\rho}{\bm}
C_1C_2^{|\bm|+n+l}\frac{1}{(n+l)!} \\
&\le 
C_1C_2^{|\bm|+n} 2^{\pair{\rho}{\bm}} 
\sum_{0\le l\le \pair{\rho}{\bm}}
\binom{\pair{\rho}{\bm}}{l}  \|(\rho+\lambda)\|^{\pair{\rho}{\bm}-l}
C_2^{l} \frac{l!}{(n+l)!} \\
&\le
C_1C_2^{|\bm|+|n|}2^{\pair{\rho}{\bm}}
(\|\rho+\lambda\|+2C_2)^{\pair{\rho}{\bm}} 
\begin{cases}
1/n! & n\ge 0, \\
|n|! & n\le 0.
\end{cases} 
\end{align*}
Thus, we obtain the estimate. 
The latter part follows from (\ref{eq:hypmodifforT}). 
\end{proof}

Let $\mO_{\lambda}^{\hbar,\hbar^{-1}}(\rho)$ 
be the set of formal power series of the form 
$\sum_{\bm\ge 0,n\in \Z} A_{\bm,n}(\lambda) x^\bm\hbar^{-n}$ 
in $\C(\lambda)[\hbar,\hbar^{-1}][\![x]\!]$ satisfying
the following conditions: 

(i) $A_{\bm,n}=0$ for $-n> \pair{\rho}{\bm}$. 

(ii) There exist positive continuous functions $B(\lambda), C(\lambda)$ 
defined on the complement of a finite subset of $\C$ such that  
\[|A_{\bm,n}(\lambda)|\le B(\lambda) C(\lambda)^{|\bm|+|n|}
 \begin{cases}1/n! & n\ge 0, \\ |n|! & n\le 0.  \end{cases} 
\]
The above lemma says that each matrix element of $\tT(x,\hbar,\lambda)$ 
is contained in $\mO_{\lambda}^{\hbar,\hbar^{-1}}(\rho)$. 

\begin{lemma}
\label{lem:subring}
The set $\mO_{\lambda}^{\hbar,\hbar^{-1}}(\rho)$ is a subring of 
$\C(\lambda)[\hbar,\hbar^{-1}][\![x]\!]$. 
\end{lemma}
\begin{proof}
We must check that $\mO_{\lambda}^{\hbar,\hbar^{-1}}(\rho)$ 
is closed under product. 
Let $\sum_{\bm\ge 0,n} A_{\bm,n} x^\bm\hbar^{-n}$, 
$\sum_{\bm\ge 0,n} B_{\bm,n} x^\bm\hbar^{-n}$ 
be in $\mO_{\lambda}^{\hbar,\hbar^{-1}}(\rho)$. 
By taking the maximum if necessary, 
we can assume that these two elements are estimated 
by the same functions $B(\lambda),C(\lambda)>1$. 
Set $C_{\bm,n}=
\sum_{\bm_1+\bm_2=\bm}\sum_{n_1+n_2=n} A_{\bm_1,n_1}B_{\bm_2,n_2}$. 
When $n\ge 0$, we have 
\begin{align*}
|C_{\bm,n}| & \le 
B(\lambda)^2 C(\lambda)^{|\bm|+|n|}
\Biggl(\sum_{\substack{n_1+n_2=n,\\ n_i\ge 0}} \frac{1}{n_1!n_2!}+
 2\sum_{i=1}^{\pair{\rho}{\bm}}C(\lambda)^{2i} \frac{i!}{(n+i)!}
\Biggr)
\sum_{\bm_1+\bm_2=\bm} 1 \\
& \le 
B(\lambda)^2 C(\lambda)^{|\bm|+|n|}
\frac{1}{n!}\left(
2^n+2\frac{C(\lambda)^{2\pair{\rho}{\bm}+2}}{C(\lambda)^2-1}\right){C'}^{|\bm|}
\end{align*}
for some $C'>0$. 
When $n\le 0$, we have 
\begin{align*}
|C_{\bm,n}|& \le 
B(\lambda)^2 C(\lambda)^{|\bm|+|n|}
\Biggl(\sum_{\substack{n_1+n_2=n, \\ n_i\le 0}} |n_1|!|n_2|! 
+2\sum_{i=1}^{\pair{\rho}{\bm}-|n|}C(\lambda)^{2i}\frac{(i+|n|)!}{i!}
\Biggr) \sum_{\bm_1+\bm_2=\bm} 1\\
& \le
B(\lambda)^2 C(\lambda)^{|\bm|+|n|}|n|!
\left(\pair{\rho}{\bm}+1+2 (1+C(\lambda)^2)^{\pair{\rho}{\bm}}\right) 
{C'}^{|\bm|}
\end{align*}
for some $C'>0$. 
Therefore, we have the desired estimate for the product. 
\end{proof}

Let $\Omega_j(x,\hbar,\lambda)$ be the connection matrix of $\tnabla$ 
given in (\ref{eq:tnabla}). 
Let $\pairr{\cdot}{\cdot}$ be the bilinear form defined in 
Proposition \ref{prop:regularity} (ii). 
\begin{lemma}
\label{lem:estofpairing}
The elements of the form 
$\pairr{\alpha}{\beta},\pairr{\alpha}{\Omega_j\beta}$
are contained in both $\C[\hbar,\lambda][\![x]\!]$ and  
$\mO_{\lambda}^{\hbar,\hbar^{-1}}(\rho)$ 
for $\alpha,\beta\in H^*(X)$.  
\end{lemma}
\begin{proof}
By Proposition \ref{prop:regularity}, we can see that 
$\pairr{\alpha}{\beta}$ and $\pairr{\alpha}{\Omega_j\beta}$ 
are in $\C[\hbar,\lambda][\![x]\!]$. 
Since we have 
\[
\pairr{\alpha}{\beta}=
\pairL{\tT(x,-\hbar,\lambda)\alpha}{\tT(x,\hbar,\lambda)\beta}_{S^1}, 
\]
we can see that $\pairr{\alpha}{\beta}$ is in 
$\mO_{\lambda}^{\hbar,\hbar^{-1}}(\rho)$ 
by Lemma \ref{lem:estofmodification}. 
On the other hand, we have 
\begin{align*}
\pairr{\alpha}{\Omega_j\beta}
&=\pairL{\tL(\tau,-\hbar)\alpha}
 {\hbar\parfrac{}{t_j}\tL(\tau,-\hbar)\beta}_{S^1} \\
&=\pairL{\tT(x,-\hbar,\lambda)\alpha}
{\nabla^{\rm cl}_j \tT(x,\hbar,\lambda)\beta}_{S^1},  
\end{align*}
where $\nabla^{\rm cl}_j=\hbar \partial_{t_j}+p_j$. 
Therefore, $\pairr{\alpha}{\Omega_j\beta}$ is also 
in $\mO_{\lambda}^{\hbar,\hbar^{-1}}(\rho)$ 
by Lemma \ref{lem:estofmodification}. 
\end{proof}
Define another ring $\mO_{\lambda}^\hbar(\rho)$ by  
\[
\mO_{\lambda}^\hbar(\rho)=
\mO_{\lambda}^{\hbar,\hbar^{-1}}(\rho)\cap \C(\lambda)[\hbar][\![x]\!].
\]
By the above lemma, $\pairr{\alpha}{\beta}$, $\pairr{\alpha}{\Omega_j\beta}$
are contained in $\mO_{\lambda}^\hbar(\rho)$. 
\begin{lemma}
\label{lem:localring}
The ring $\mO_{\lambda}^\hbar(\rho)$ is a local ring. 
\end{lemma}
\begin{proof}
Let $\sum_{\bm\ge 0}\sum_{n=0}^{\pair{\rho}{\bm}} 
A_{\bm,n}(\lambda)x^\bm \hbar^n$ be in $\mO_{\lambda}^\hbar(\rho)$. 
We will show that if $A_{0,0}(\lambda)=1$, 
this element is invertible in $\mO_{\lambda}^\hbar(\rho)$. 
We can assume that $|A_{\bm,n}|\le B(\lambda)C(\lambda)^{|\bm|+n}n!$ 
for some positive functions $B(\lambda),C(\lambda)>1$. 
Set $1+\sum_{\bm>0,n\ge 0}B_{\bm,n}x^\bm \hbar^n=
(1+\sum_{\bm>0,n\ge 0}A_{\bm,n}x^\bm\hbar^n)^{-1}$. 
Then we have 
\begin{align*}
|B_{\bm,n}|& =
\left|\sum_{l\ge 1}(-1)^l \sum_{\substack{\bm=\bm_1+\cdots+\bm_l \\ \bm_i>0}}
\sum_{\substack{n=n_1+\cdots +n_l\\ 0\le n_i\le \pair{\rho}{\bm_i}}}
A_{\bm_1,n_1}\cdots A_{\bm_l,n_l}\right| \\
&\le \sum_{l\ge 1}\sum_{\substack{\bm=\bm_1+\cdots+\bm_l \\ \bm_i>0}}
\sum_{\substack{n=n_1+\cdots +n_l\\ 0\le n_i\le \pair{\rho}{\bm_i}}}
B(\lambda)^l C(\lambda)^{|\bm|+n}n_1!\cdots n_l! \\
&\le B(\lambda)^{|\bm|} C(\lambda)^{|\bm|+n}n!
\sum_{l\ge 1}\sum_{\substack{\bm=\bm_1+\cdots+\bm_l \\ \bm_i>0}}
\binom{n+l-1}{n} \\
&\le B(\lambda)^{|\bm|}C(\lambda)^{|\bm|+n}2^{n+|\bm|-1} n! 
\sum_{l\ge 1}\sum_{\substack{\bm=\bm_1+\cdots+\bm_l \\ \bm_i>0}} 1,
\end{align*}
where we used $l\le |\bm|$. 
The summation factor is of exponential order in $|\bm|$ because 
\[
1+\sum_{l\ge 1}\sum_{\substack{\bm=\bm_1+\dots+\bm_l,\\ \bm_i>0}} y^{|\bm|}
=\frac{1}{1-\sum_{\bm>0}y^{|\bm|}}
=\frac{1}{1-(1/(1-y)^s-1)}
\]
is analytic around $y=0$. 
\end{proof}
\begin{proposition}
Each entry of the connection matrix $\Omega_j(x,\hbar,\lambda)$ 
is contained in the ring $\mO_{\lambda}^\hbar(\rho)$. 
\end{proposition}
\begin{proof}
By Lemma \ref{lem:estofpairing}, $\eta_{kl}=\pairr{p_k}{p_l}$ 
and $\pairr{p_k}{\Omega_j p_l}$ belong to $\mO_{\lambda}^\hbar(\rho)$. 
It is easy to see that 
$\lim_{q\to 0}\eta_{kl}=\pairL{p_k}{p_l}_{S^1}$. 
Because the matrix $(\pairL{p_k}{p_l}_{S^1})$ is invertible  
in $\Mat(s+1,\C(\lambda))$ 
and $\mO_{\lambda}^{\hbar}(\rho)$ is a local ring,  
the matrix $(\eta_{kl})$ is invertible in 
$\Mat(s+1,\mO_{\lambda}^{\hbar}(\rho))$. 
Therefore, matrix elements 
$\Omega_{j;kl}=\sum_i\eta^{ki}\pairr{p_i}{\Omega_jp_l}$ 
are in $\mO_{\lambda}^\hbar(\rho)$. 
\end{proof}

\subsection{Gauge fixing}
In this subsection, we will find a gauge transformation $g_2$ 
which changes the connection matrices into $\hbar$-independent ones. 

Let $\mO_{\lambda}^\hbar$ be the set of formal power series 
$\sum_{\bm,n\ge 0}A_{\bm,n}(\lambda)x^\bm\hbar^n$ 
in $\C(\lambda)[\hbar][\![x]\!]$ satisfying the following conditions: 

(i) $A_{\bm,n}=0$ if  $|\bm|_0=0,n>0$, 

(ii) There exists positive continuous functions 
$B(\lambda)$ and $C(\lambda)$ defined on the complement of 
a finite subset of $\C$ such that 
\[
|A_{\bm,n}(\lambda)|\le B(\lambda)C(\lambda)^{|\bm|+n}|\bm|_0^n, 
\]
where $|\bm|_0=\sum_{a=1}^r m_a (\le |\bm|)$. 
It is easy to see that $\mO_{\lambda}^\hbar(\rho)$ 
is contained in  $\mO_{\lambda}^\hbar$.
\begin{lemma}
$\mO_{\lambda}^\hbar$ is a local ring. 
\end{lemma}
We omit the proof because 
it is similar to Lemma \ref{lem:localring}.

\begin{proposition}
\label{prop:gaugefixing}
There exists a unique gauge transformation 
$g_2$ with entries in $\mO_{\lambda}^{\hbar}$ 
such that $g_2|_{q=0}=\id$ and that 
the new connection matrix $\hOmega_j$ 
of $g_2^*\tnabla=\hbar d+ \sum_{j=0}^s\hOmega_j dt_j$
is $\hbar$-independent. 
Moreover, $\hOmega_j$ is convergent. 
\end{proposition}
This proposition is considered to be a general gauge fixing lemma. 
It is applicable to any flat connection of the form 
$\hbar d+\Omega$  which is defined over $\mO_{\lambda}^\hbar$,  
regular singular along $q_1q_2\cdots q_r=0$ and 
whose residue matrices at $q=0$ are nilpotent. 
In \cite{iritani-genmir}, Theorem 4.6, 
we showed the existence and the uniqueness of $g_2$ 
in $\End(H^*(X))\otimes\C[\hbar,\lambda][\![x]\!]$
by using a formal Birkhoff factorization.
Here, we will show that $g_2$ also belongs to 
$\End(H^*(X))\otimes\mO_\lambda^\hbar$. 
This gauge fixing can be considered as a procedure of 
{\it renormalization}. 
A divergent connection can be renormalized by $g_2$  
to yield a finite (convergent) result. 
\begin{proof}
Once we establish the existence of 
$g_2$ in $\End(H^*(X))\otimes\mO_{\lambda}^\hbar$, 
we can see that $\hOmega_j$ is convergent because 
it is contained in  $\End(H^*(X))\otimes\mO_{\lambda}^\hbar$ 
and $\hbar$-independent at the same time.  
Thus, it suffices to solve for $g_2$ in 
$\End(H^*(X))\otimes\mO_{\lambda}^\hbar$.  
Set $\Omega_j=\sum_{\bd,n\ge 0}
\Omega_{j;\bd,n}(t_{r+1},\cdots,t_s,\lambda)q^\bd \hbar^n$. 
Because $\Omega_j$ is defined over $\mO_{\lambda}^{\hbar}$, 
there exist a neighborhood $U\subset \C^{s-r}$ of $0$ 
and continuous functions $B(\lambda),C(\lambda)>0$ 
such that 
$\|\Omega_{j;\bd,n}(t,\lambda)\|
\le B(\lambda)C(\lambda)^{|\bd|+n}|\bd|^n$
for $(t_{r+1},\cdots,t_s)\in U$ and $0\le j\le s$. 
From now on, we omit $\lambda$ and $(t_{r+1},\cdots, t_s)$ in the notation, 
but $(t_{r+1},\cdots,t_s)$ is always assumed 
to be in $U$.  
We set 
$g_2=\sum_{\bd,n\ge 0} G_{\bd,n}q^\bd\hbar^n,
\hOmega_j=\sum_{\bd,n\ge 0}\hOmega_{j;\bd}q^\bd$.
By expanding the relation 
\[ 
g_2 \hOmega_a=\Omega_a g_2 +\hbar q_a\parfrac{g_2}{q_a},\quad 
1\le a\le r,  
\]
we obtain the following equations:
\begin{align}
\label{eq:inductionnge1}
d_a G_{\bd,n-1}+\ad(p_a) G_{\bd,n}
&=\sum_{\substack{\bd_1+\bd_2=\bd, \\ \bd_1>0,\bd_2>0 }} 
  G_{\bd_2,n}\hOmega_{a;\bd_1}
  -\sum_{\substack{\bd_1+\bd_2=\bd,\\ n=n_1+n_2 \\\bd_1>0}}
  \Omega_{a;\bd_1,n_1} G_{\bd_2,n_2}\quad  (n\ge 1), \\
\label{eq:inductionn=0}
\ad(p_a)G_{\bd,0}
&=\hOmega_{a;\bd}+
  \sum_{\substack{\bd_1+\bd_2=\bd, \\ \bd_1>0,\bd_2>0 }}
  G_{\bd_2,0}\hOmega_{a;\bd_1}
  -\sum_{\substack{\bd_1+\bd_2=\bd,\\\bd_1>0}}
  \Omega_{a;\bd_1,0}G_{\bd_2,0},  
\end{align}
where we used $G_{0,n}=\delta_{0,n}$, $\hOmega_{a;0}=p_a$ and 
$\Omega_{a;0,n}=\delta_{0,n}p_a$. 
Note that $\Omega_{a;\bd,n}$ is known and 
$G_{\bd,n}$ and $\hOmega_{a;\bd}$ are unknown. 
Assume by induction that we know 
$G_{\bd',n}$ and $\hOmega_{a;\bd'}$ for all $\bd'$ with $|\bd'|<\ovd$. 
For a multi-index $\bd$ with $|\bd|=\ovd$, 
we first solve for $G_{\bd,n}$ for all $n\ge 0$ by using 
(\ref{eq:inductionnge1}), and then we 
solve for $\hOmega_{a;\bd}$ by using (\ref{eq:inductionn=0}). 

More precisely, we must solve for $G_{\bd,n}, \hOmega_{a;\bd}$
with estimates. 
Introduce the following notation: 
\begin{gather*}
g_{d,n}=\frac{1}{d^n}\sum_{|\bd|=d}\|G_{\bd,n}\|,\quad 
\omega_{a;d,n}=\frac{1}{d^n}\sum_{|\bd|=d}\|\Omega_{a;\bd,n}\|,\quad 
\homega_{a;d}=\sum_{|\bd|=d}\|\hOmega_{a;\bd}\|, \\
H_{a;\bd,n}=
\text{the right hand side of (\ref{eq:inductionnge1})}, \quad 
h_{d,n}=\max_{1\le a\le r}
\frac{1}{d^n}\sum_{|\bd|=d}\|H_{a;\bd,n}\|\quad (n\ge 1).
\end{gather*}
There exist positive $A_1,B_1$ such that 
$\omega_{a;d,n}\le A_1B_1^{d+n}$. 
Assume by induction that 
\begin{equation}
\label{eq:assumeest}
g_{d,n}\le \frac{B_2^{d}}{(d+1)^M}B_3^n,\quad 
\homega_{a;d}\le A_2\frac{B_2^d}{(d+1)^M}\quad (1\le a\le r) 
\end{equation}
hold for all $d<\ovd$. 
We take $A_2$ so that $\homega_{a;0}=\|p_a\|\le A_2$ holds. 
Then (\ref{eq:assumeest}) is valid for $d=0$ 
because $g_{0,n}=\delta_{0,n}$. 
We will specify $B_2,B_3,M$ later. 
Take a $\bd$ with $|\bd|=\ovd$.  
Let $a(\bd)$ be an index such that $d_{a(\bd)}=\max\{d_1,\dots,d_r\}$. 
Let $C$ be a constant satisfying $\|\ad(p_a)\|\le C$ 
for $1\le a\le r$. 
By (\ref{eq:inductionnge1}), we have 
\begin{align*}
G_{\bd,n}
&=\frac{1}{d_{a(\bd)}}(H_{a(\bd);\bd,n+1}-\ad(p_{a(\bd)})G_{\bd,n+1}) \\
&=\frac{1}{d_{a(\bd)}}H_{a(\bd);\bd,n+1}
  -\frac{\ad(p_{a(\bd)})}{d_{a(\bd)}^2}H_{a(\bd);\bd,n+2}+
  \cdots +\frac{\ad(p_{a(\bd)})^{2N}}{d_{a(\bd)}^{2N+1}}
  H_{a(\bd);\bd,n+2N+1},
\end{align*}
where $N=\dim_{\C}X$ and 
we used $\ad(p_{a(\bd)})^{2N+1}=0$. 
By using $(|\bd|/d_{a(\bd)})\le r$, we have 
\begin{equation}
\label{eq:gdnest}
g_{\ovd,n}\le r h_{\ovd,n+1}+r^2C h_{\ovd,n+2}+\cdots+
r^{2N+1}C^{2N}h_{\ovd,n+2N+1}.
\end{equation}
Also we have 
\begin{align}
\label{eq:hdnest}
h_{\ovd,n}&\le \max_{1\le a\le r}\left(
\sum_{\substack{|\bd_1|+|\bd_2|=\ovd,\\ \bd_1>0,\bd_2>0}}
\frac{1}{\ovd^n}\|G_{\bd_2,n}\|\|\hOmega_{a;\bd_1}\| 
+\sum_{\substack{|\bd_1|+|\bd_2|=\ovd,\\ n_1+n_2=n, \bd_1>0}}
\frac{1}{\ovd^{n_1}}\|\Omega_{a;\bd_1,n_1}\|
\frac{1}{\ovd^{n_2}}\|G_{\bd_2,n_2}\|\right) \\
\nonumber 
&\le \max_{1\le a\le r}\left(
\sum_{d_1=1}^{\ovd-1}g_{\ovd-d_1,n}\homega_{a;d_1}
+\sum_{d_1=1}^{\ovd} \sum_{k=0}^n \omega_{a;d_1,k}g_{\ovd-d_1,n-k}\right).
\end{align}
By using (\ref{eq:inductionn=0}), we have 
\begin{align}
\label{eq:homegaadest}
\homega_{a;\ovd}
&\le C \sum_{|\bd|=\ovd}\|G_{\bd,0}\| 
+\sum_{\substack{|\bd_1|+|\bd_2|=\ovd,\\ \bd_1>0,\bd_2>0}}
\|G_{\bd_2,0}\|\|\hOmega_{a;\bd_1}\| 
+\sum_{\substack{|\bd_1|+|\bd_2|=\ovd,\\ \bd_1>0}}
\|\Omega_{a;\bd_1,0}\|\|G_{\bd_2,0}\|  \\
\nonumber 
&\le C g_{\ovd,0} 
+\sum_{d_1=1}^{\ovd-1} g_{\ovd-d_1,0}\homega_{a;d_1} 
+\sum_{d_1=1}^{\ovd} \omega_{a;d_1,0} g_{\ovd-d_1,0}.
\end{align}
By (\ref{eq:hdnest}) and the assumption (\ref{eq:assumeest}), we obtain 
\begin{align}
\label{eq:hdnest2}
h_{\ovd,n}& \le A_2B_2^{\ovd}B_3^n \sum_{i=1}^{\ovd-1}
\frac{1}{(\ovd-i+1)^M(i+1)^M} 
+A_1B_2^{\ovd}B_3^n 
\sum_{i=1}^\ovd\sum_{k=0}^n \left(\frac{B_1}{B_2}\right)^i 
\left(\frac{B_1}{B_3}\right)^k\frac{1}{(\ovd-i+1)^M} \\
\nonumber 
& \le \left(A_2\varepsilon_1(M)
+\frac{A_1}{1-B_1/B_3} \varepsilon_2(B_2,M)\right)
\frac{B_2^{\ovd}}{(\ovd+1)^M} B_3^n,  
\end{align}
where we set 
\[
\varepsilon_1(M)=\sup_{\ovd\ge 1}\sum_{i=1}^{\ovd-1}
\frac{(\ovd+1)^M}{(\ovd-i+1)^M(i+1)^M},\quad 
\varepsilon_2(B_2,M)=\sup_{\ovd\ge 1}\sum_{i=1}^\ovd 
\left(\frac{B_1}{B_2}\right)^i \frac{(\ovd+1)^M}{(\ovd-i+1)^M}.
\]
By (\ref{eq:gdnest}) and (\ref{eq:hdnest2}), we obtain 
\begin{equation}
\label{eq:gdnest2}
g_{\ovd,n}\le  
\varepsilon_3(B_2,B_3,M) \frac{B_2^{\ovd}}{(\ovd+1)^M} B_3^n,
\end{equation} 
where we set 
\[
\varepsilon_3(B_2,B_3,M)= 
\left(A_2\varepsilon_1(M)
+\frac{A_1}{1-B_1/B_3} \varepsilon_2(B_2,M)\right)
rB_3\frac{(rB_3C)^{2N+1}-1}{rB_3C-1}.  
\]
By (\ref{eq:gdnest2}) and (\ref{eq:homegaadest}), we have 
\begin{equation}
\label{eq:homegaadest2}
\homega_{a;\ovd}\le 
\varepsilon_3(B_2,B_3,M)(C+A_2\varepsilon_1(M)+A_1\varepsilon_2(B_2,M))
\frac{B_2^{\ovd}}{(\ovd+1)^M}.
\end{equation}
In order to complete the induction step, we need to specify 
the parameters $B_2,B_3,M$.  
First we set $B_3=2B_1$. 
To choose $B_2$ and $M$, we need the following lemma: 
\begin{lemma}
\label{lem:epsilongoestozero}
$\lim_{M\to \infty}\varepsilon_1(M)=0,\ 
\lim_{B_2 \to \infty}\varepsilon_2(B_2,M)=0$.
\end{lemma}
The proof will be given in the Appendix. 
For sufficiently large $M$, we have 
\[
A_2\varepsilon_1(M) rB_3\frac{(rB_3C)^{2N+1}-1}{rB_3C-1} <
\min\left\{\frac{1}{2},\frac{A_2}{3C}\right\} \text{ and } 
\varepsilon_1(M)<\frac{1}{3}. 
\]
Next, for sufficiently large $B_2$, we have  
\[
\frac{A_1}{1-B_1/B_3}
\varepsilon_2(B_2,M) rB_3\frac{(rB_3C)^{2N+1}-1}{rB_3C-1}<\frac{1}{2}
\text{ and }A_1 \varepsilon_2(B_2,M) <\frac{A_2}{3}.
\]
Now, it is easy to check that $\varepsilon_3(B_2,B_3,M)<1$ 
and $\varepsilon_3(B_2,B_3,M)
(C+A_2\varepsilon_1(M)+A_1\varepsilon_2(B_2,M))<A_2$. 
Therefore, by (\ref{eq:gdnest2}) and (\ref{eq:homegaadest2}), 
we complete the induction step. 
\end{proof}

\begin{proof}[Proof of Theorem \ref{thm:mainthm}] 
As explained at the end of the section 
\ref{sect:quantumLefschetz}, 
in order to obtain the structure constants of $QH^*_{S^1}(X,\mL)$ 
from $g_2^*\tnabla$, 
it suffices to find a new coordinate system $\hx=
\{\hatt_0,\hq_1=\exp(\hatt_1),\dots,\hq_r=\exp(\hatt_r),
\hatt_{r+1},\dots,\hatt_s\}$ 
such that the connection matrix 
$\hOmega_{\hat{j}}$ defined by 
$g_2^*\tnabla=\hbar d+\sum_{j=0}^s \hOmega_j dt_j
=\hbar d+\sum_{j=0}\hOmega_{\hat{j}}d\hatt_j$
satisfies $\hOmega_{\hat{j}}1=p_j$. 
Then $\hOmega_{\hat{j}}$ gives the twisted quantum product $p_j*_{\mL}$. 
Because $\hOmega_j$ is already convergent, 
new coordinates $\hx$ also become convergent functions in $x$. 
\end{proof}

\section{Mirror symmetry for non-nef toric varieties}
\label{sect:mirrorsymmetryfornonneftoricvarieties} 
In this section, we will study mirror symmetry and 
the quantum cohomology of a not necessarily nef toric variety. 
In \cite{iritani-genmir}, 
we showed that the quantum $D$-module of $QDM^*(X)$ of a toric variety $X$ 
can be reconstructed from the equivariant Floer cohomology 
$\FH$ \cite{givental-homologicalgeom,iritani-EFC} 
by a generalized mirror transformation. 
There, the quantum $D$-module $QDM^*(X)$ was described in terms of 
hypergeometric series ($I$-function). 
Here, we describe $QDM^*(X)$ in terms of 
the mirror oscillatory integral introduced by Givental   
\cite{givental-ICM,givental-mirrorthm-toric}. 
The oscillatory integral satisfies (a generalized version of) 
the Mellin system of hypergeometric differential equations.  
We will show that the equivariant Floer cohomology $\FH$ is isomorphic 
to the {\it $q$-adic completion} of the Mellin system. 
Then, results in \cite{iritani-genmir} can be restated 
as $QDM^*(X)$ restricted to some non-linear subspace of $H^*(X)$ 
is isomorphic to the completion of the Mellin system. 
Using a method similar to section \ref{sect:proofofmainthm}, 
we will show the convergence of 
the quantum cohomology of toric varieties. 
By using mirror symmetry, we will also show the semisimplicity. 

\subsection{Mirrors and the Mellin system}
\label{sect:mirrorsandtheMellinsystem}
Let $X$ be a smooth projective toric variety 
defined by a fan $\Sigma\subset \R^N$.  
Take a nef integral basis $\{p_1,\dots,p_r\}$ of $H^2(X,\Z)$.
Let $D_1,\dots,D_{r+N}$ be all the torus invariant prime divisors 
and $\td_1,\dots,\td_{r+N}\in H^2(X,\Z)$ be 
their Poincar\'{e} duals.  
We write $w_i=\sum_{a=1}^r m_{ia}p_a$. 
We can recover $X$ from the data $w_i$ and K\"{a}hler class 
$\eta_X$ in $H^2(X,\R)$. 
Set $\mathcal{B}=\{I\subset\{1,\dots,N+r\}\;|\;\eta_X\notin 
\sum_{i\in I}\R_{>0}\td_i\}$. 
Then we have 
\[
X=\C^{r+N}_{\mathcal{B}}/(\C^*)^r,\quad 
\C^{r+N}_{\mathcal{B}}:=\C^{r+N}\setminus\bigcup_{I\in\mathcal{B}}
  \{(z_1,\dots,z_{r+N})\;|\; z_i=0 \text{ for }i\notin I\}.
\] 
Here, $(\C^*)^r$ acts on $\C^{r+N}$ as 
$(z_1,\dots,z_{r+N})\mapsto (t^{\td_1}z_1,\dots, t^{\td_{r+N}}z_{r+N})$, 
where $t^{\td_i}=\prod_{a=1}^r t_a^{m_{ia}}$. 
The divisor $D_i$ corresponds to $\{z_i=0\}$. 
Let $\mc_i$ denotes a coordinate of the mirror corresponding 
to the class $\td_i$. 
Let $\pi\colon Y=(\C^*)^{r+N}\rightarrow (\C^*)^r$ 
be a family of algebraic tori defined by 
$\pi(\mc_1,\dots,\mc_{r+N})=(q_1,\dots,q_r)$, 
$q_a=\prod_{i=1}^{r+N} \mc_i^{m_{ia}}$. 
Define a function $F(\mc)$ on $Y$ as  $F(\mc)=\mc_1+\cdots+\mc_{r+N}$.
We write $Y_{q}=\pi^{-1}(q)$ and $F_q=F|_{Y_q}$ for $q\in (\C^*)^r$. 
The mirror oscillatory integral is given by 
\[
\mI_{\Gamma}(q,\hbar)=\int_{\Gamma_q\subset Y_q} e^{F_q/\hbar}\omega_q,  
\]
where $\Gamma_q$ is a non-compact real $N$-cycle in $Y_q$ 
such that the integral converges.  
The holomorphic volume form $\omega_q$ is given below. 
Take coordinates $(s_1,\cdots,s_N)$  of fibers $Y_q$ of the form  
$s_b=\prod_{i=1}^{r+N}\mc_i^{m'_{ib}}$, 
where the $(r+N)\times N$ matrix 
$(m'_{ib})\colon \Z^{r+N}\rightarrow \Z^N$ gives a splitting of 
the exact sequence 
\begin{equation}
\label{eq:toricexactseq}
\begin{CD}
0 @>>> \Z^N @>\ker (m_{ia}) >> \Z^{r+N} @> (m_{ia}) >> 
\Z^r\cong H^2(X,\Z) @>>> 0. 
\end{CD}
\end{equation}
By this splitting, we can write 
$\mc_i=\prod_{a=1}^r q_a^{l_{ai}}\prod_{b=1}^N s_b^{x_{ib}}$. 
Here, $\vx_i=(x_{i1},\dots,x_{iN}) \in \Z^N$ gives the primitive generator  
of the $i$-th one dimensional cone of the fan $\Sigma$. 
The matrix  $(l_{ai})\colon \Z^r\rightarrow \Z^{r+N}$ 
also gives a splitting of the exact sequence (\ref{eq:toricexactseq}). 
We can choose this splitting so that $l_{ai}\ge 0$ because 
$\{p_1,\dots,p_r\}$ is a nef basis.  
We define $\omega_q$ by  
\[
\omega_q=
\frac{d\log \mc_1\wedge\cdots\wedge d\log \mc_{r+N}}
     {d\log q_1\wedge\cdots\wedge d\log q_r} 
=d \log s_1\wedge \cdots \wedge d\log s_N|_{Y_q}. 
\] 
This is independent of a choice of the splitting. 
\begin{proposition}
\label{prop:noneqmirrordiffeq}
The oscillatory integral $\mI_{\Gamma}(q,\hbar)$ 
satisfies  $\mP_\bd \mI_{\Gamma}(q,\hbar)=0$ for all $\bd\in \Z^r$, 
where $\mP_\bd$ is a differential operator defined by
\[
\mP_{\bd}=
q^\bd \prod_{\pair{\td_i}{\bd}<0} \prod_{k=0}^{-\pair{\td_i}{\bd}-1}
\left(\sum_{a=1}^r m_{ia} \hbar\partial_a - k\hbar\right)
-\prod_{\pair{\td_i}{\bd}>0} 
\prod_{k=0}^{\pair{\td_i}{\bd}-1}
\left(\sum_{a=1}^r m_{ia} \hbar\partial_a - k\hbar\right), 
\]
where $\partial_a=q_a\partial/\partial q_a$. 
\end{proposition} 
This proposition may be well-known, 
but we include a proof for completeness.  
\begin{proof}
First we have 
\begin{align*}
\int_{\Gamma} \hbar\parfrac{e^{F_q/\hbar}}{\log \mc_i}\omega_q
&= \int_{\Gamma} 
\left(\sum_{a=1}^r m_{ia} \hbar\partial_a e^{F_q/\hbar}+
 \sum_{b=1}^N m'_{ib} \hbar\parfrac{e^{F_q/\hbar}}{\log s_b}\right)
 \prod_{b} d \log s_b  \\
&= \left(\sum_{a=1}^r m_{ia} \hbar\partial_a\right) 
   \int_{\Gamma} e^{F_q/\hbar}\prod_{b} d \log s_b  
  + \sum_{c=1}^N m'_{ic} \int_{\Gamma} \hbar d \left({e^{F_q/\hbar}}
     \prod_{b\neq c} d \log s_b\right)  \\
&= \left(\sum_{a=1}^r m_{ia} \hbar\partial_a\right) 
   \int_{\Gamma} e^{F_q/\hbar}\omega_q. 
\end{align*}
By using this, we have 
\begin{align*}
\mP_\bd \mI(q,\hbar) &= \int_{\Gamma} 
\Bigl(q^\bd \prod_{\pair{\td_i}{\bd}<0} \prod_{k=0}^{-\pair{\td_i}{\bd}-1}
 (\hbar \mc_i\parfrac{}{\mc_i}-k\hbar)-
\prod_{\pair{\td_i}{\bd}>0} \prod_{k=0}^{\pair{\td_i}{\bd}-1}
 (\hbar \mc_i\parfrac{}{\mc_i}-k\hbar)\Bigr) e^{F_q/\hbar}\omega_q \\
&=\int_{\Gamma} 
\Bigl(q^\bd \prod_{\pair{\td_i}{\bd}<0} \mc_i^{-\pair{\td_i}{\bd}}-
\prod_{\pair{\td_i}{\bd}>0} \mc_i^{\pair{\td_i}{\bd}}\Bigr)
e^{F_q/\hbar}\omega_q=0,  
\end{align*}
where we used $q^\bd=\prod_{i}\mc_i^{\pair{\td_i}{\bd}}$. 
\end{proof}
We call the system of hypergeometric differential equations 
$\{\mP_\bd \mI(q,\hbar)=0\;|\; \bd\in \Z^r\}$ the Mellin system. 
Let $M_{\rm Mell}$ be the $D$-module corresponding to the Mellin system 
\[
M_{\rm Mell}=
\C\langle q^{\pm},\hbar\partial, \hbar\rangle \big/I_{\rm Mell},
\quad  
I_{\rm Mell}=
\sum_{\bd\in \Z^r} \C\langle q^{\pm},\hbar\partial,\hbar\rangle\mP_\bd,  
\]
where $q^{\pm}$ and $\hbar\partial$ 
are shorthand for $q_1^{\pm},\dots,q_r^{\pm}$ and 
$\hbar q_1\partial/\partial q_1,\dots,\hbar q_r\partial/\partial q_r$. 
The oscillatory integral $\mI_{\Gamma}(q,\hbar)$ gives 
a solution of $M_{\rm Mell}$ for each non-compact cycle $\Gamma$. 

\subsection{$S^1$-equivariant Floer cohomology}
\label{sect:S^1-equivariantFloercohomology}
We review the algebraic construction of the equivariant 
Floer cohomology $\FH$ for a toric variety $X$ briefly
(see \cite{iritani-EFC} for detail). 
For each $\bd\in \Z^r$, we put 
$H^*_{S^1}(L_\bd^\infty)=\C[P_1,\dots,P_r,\hbar]$. 
This is an algebra over $H^*_{S^1}({\rm pt})=\C[\hbar]$. 
When $\pair{\td_i}{\bd}\ge \pair{\td_i}{\bd'}$ for all $i$, 
we define a push-forward map $i_{\bd,\bd'}\colon 
H^*_{S^1}(L_\bd^\infty)\rightarrow H^*_{S^1}(L_{\bd'}^\infty)$ 
by $i_{\bd,\bd'}(\alpha)=\alpha \prod_{i=1}^N 
\prod_{k=\pair{\td_i}{\bd'}}^{\pair{\td_i}{\bd}-1}\tdD_{i,k}$, 
where $\tdD_{i,k}=\sum_{a=1}^r m_{ia}P_a -k\hbar$.  
Then, we have an inductive system 
$(H^*_{S^1}(L_\bd^\infty),i_{\bd,\bd'})$. 
Let $H^{\infty/2}_{S^1}=\varinjlim_{\bd} H^*_{S^1}(L_{\bd}^\infty)$ 
be its direct limit.  
We define an $H^*_{S^1}({\rm pt})$-algebra homomorphism  
$Q^{\bd}=\prod_{a=1}^r Q_a^{d_a}
\colon H^*_{S^1}(L_{\bd'}^\infty)\rightarrow H^*_{S^1}(L_{\bd'+\bd}^\infty)$ 
by $Q^{\bd}(P_a)=P_a-d_a\hbar$. 
This is compatible with the direct limit and 
we have a module homomorphism 
$Q^\bd\colon H^{\infty/2}_{S^1}\rightarrow H^{\infty/2}_{S^1}$. 
We can check that the multiplication by $P_a$ and the action of $Q_b$ 
on $H^{\infty/2}_{S^1}$ satisfies 
the commutation relation $[P_a,Q_b]=\hbar\delta_{ab}Q_b$. 
Hence, $H^{\infty/2}_{S^1}$ has a $D$-module structure 
when we regard $P_a$ as a differential operator 
$\hbar Q_a(\partial/\partial Q_a)$. 
Let $\Delta_\bd\in H_{S^1}^{\infty/2}$ be the image of 
$1\in H^*_{S^1}(L_{\bd}^\infty)$. 
We also write $\Delta=\Delta_0$. 
Let $FH_0$ be the submodule 
$\C\langle P_1,\dots,P_r,Q_1,\dots,Q_r,\hbar\rangle \Delta$ 
of $H_{S^1}^{\infty/2}$. 
This $FH_0$ has a natural $Q$-adic topology and we define 
$\FH$ as the completion of $FH_0$: 
\[
FH^*_{S^1}= \widehat{FH_0} =\varprojlim_{n} FH_0/\frm^n FH_0,
\]
where $\frm=\sum_{a}Q_a\C[Q,\hbar]$. 
This becomes a module over 
$\C[\hbar][\![Q_1,\dots,Q_r]\!]\langle P_1,\dots,P_r\rangle$. 

\begin{proposition}
\label{prop:H^infty/2=Mell}
The $D$-module $H_{S^1}^{\infty/2}$ is generated by $\Delta$ 
as a $\C\langle Q^{\pm},P,\hbar\rangle$-module and 
all the relations are generated by $\widetilde{\mP}_{\bd}
=Q^\bd\prod_{\pair{\td_i}{\bd}<0}\prod_{k=0}^{-\pair{\td_i}{\bd}-1} 
\tdD_{i,k}-\prod_{\pair{\td_i}{\bd}>0}\prod_{k=0}^{\pair{\td_i}{\bd}-1}
\tdD_{i,k}$
for $\bd\in \Z^r$. 
Therefore, $H^{\infty/2}_{S^1}$ is isomorphic to $M_{\rm Mell}$. 
\end{proposition} 
\begin{proof}
By the relation $\Delta_\bd=Q^\bd\Delta$, 
$H_{S^1}^{\infty/2}$ is generated by $\Delta$. 
It is easy to check that $\widetilde{\mP}_\bd\Delta=0$. 
Assume that $f(P,Q,\hbar)\Delta=0$. 
We set $f(P,Q,\hbar)=\sum_{i}f_i(P,\hbar)Q^{\bd_i}$. 
There exists $\bd$ such that $\pair{\td_j}{\bd_i+\bd}>0$ for all $i,j$. 
Then we have 
\begin{align}
\label{eq:relationofsemi-infcoh}
Q^{\bd}f(P,Q,\hbar)&=
\sum_{i}f_i(P_a-\pair{p_a}{\bd}\hbar,\hbar)Q^{\bd_i+\bd}\\
\nonumber 
&=\sum_{i}f_i(P_a-\pair{p_a}{\bd}\hbar,\hbar)
\Bigl( \widetilde{\mP}_{\bd_i+\bd}-
\prod_{j}\prod_{k=0}^{\pair{\td_j}{\bd_i+\bd}-1}\tdD_{j,k}\Bigr).   
\end{align}
When applying this to $\Delta$, we have 
$\sum_i f_i(P_a-\pair{p_a}{\bd}\hbar,\hbar)
\prod_{j}\prod_{k=0}^{\pair{\td_j}{\bd_i+\bd}-1}\tdD_{j,k}\Delta=0$. 
Because the canonical map $H^*_{S^1}(L_0^\infty)\rightarrow H_{S^1}^\infty$
is injective, we have 
\[
\sum_i f_i(P_a-\pair{p_a}{\bd}\hbar,\hbar)
\prod_{j}\prod_{k=0}^{\pair{\td_j}{\bd_i+\bd}-1}\tdD_{j,k}=0.
\]
Therefore, by (\ref{eq:relationofsemi-infcoh}), we have 
$f(P,Q,\hbar)=\sum_{i} f_i(P,\hbar) Q^{-\bd}\widetilde{\mP}_{\bd_i+\bd}$. 
\end{proof}
\begin{corollary}
Under the correspondence $q_a\mapsto Q_a$, 
$\hbar\partial_a\mapsto P_a$, we have
\[ 
FH_0\cong \C\langle q,\hbar\partial,\hbar\rangle/I_{\rm poly},
\quad 
I_{\rm poly}=  I_{\rm Mell}
 \cap \C\langle q,\hbar\partial,\hbar\rangle.
\]
\end{corollary}
Hereafter, we identify $Q_a$ and $P_a$ with $q_a$ and $\hbar\partial_a$. 
From this, we can describe $\FH$ as follows.  
\begin{proposition}
\label{prop:FH=completedMell}
\[
\FH=\widehat{FH}_0\cong 
\C[\hbar ][\![q]\!]\langle \hbar\partial \rangle / 
\ov{I}_{\rm poly},   
\]
where $\ov{I}_{\rm poly}\subset\C[\hbar][\![q]\!]\langle\hbar\partial\rangle$ 
is the closure of $I_{\rm poly}$ with respect to
the $q$-adic topology.  
\end{proposition}  
\begin{proof}
In \cite{iritani-EFC}, section 4.4, it is proved that $\FH$ 
is generated by $\Delta$ over $\C[\hbar][\![Q]\!]\langle P\rangle$.
Therefore we have a surjection 
$\C[\hbar][\![Q]\!]\langle P\rangle\rightarrow \FH$. 
Assume $f(P,Q,\hbar)=\sum_{\bd\ge 0}f_{\bd}(P,\hbar)Q^{\bd} 
\in \C[\hbar][\![Q]\!]\langle P\rangle$ 
satisfies $f(P,Q,\hbar)\Delta=0$.  
Set $f_n=\sum_{|\bd|\le n} f_{\bd}(P,\hbar)Q^{\bd}$. 
Then,  
\[
f_n\Delta=-\sum_{|\bd|>n} f_{\bd}(P,\hbar)Q^\bd\Delta
\]
belongs to $FH_0\cap \frm^{n+1} \widehat{FH}_0
=\frm^{n+1} FH_0$. 
Therefore, there exists $g_n$ in $\frm^{n+1}\C\langle Q,P,\hbar\rangle$
such that 
$f_n\Delta=g_n \Delta$. 
Because $f_n-g_n\in I_{\rm poly}$ and $f_n-g_n \to f$ 
as $n\to \infty$, we have 
$f\in \ov{I}_{\rm poly}$. 
\end{proof}
Later, we will see that $\ov{I}_{\rm poly}$ 
does not necessarily coincide with 
$\C[\hbar][\![q]\!]\langle \hbar\partial \rangle I_{\rm poly}$. 
In such a case, we need to add non-algebraic differential equations.
We set $\mO_{\rm small}^{\hbar}=\mO_{\lambda}^{\hbar}\cap \C[\hbar][\![q]\!]$ 
($\mO_{\lambda}^{\hbar}$ was introduced in 
section \ref{sect:proofofmainthm}): 
\[
\mO^{\hbar}_{\rm small} = 
\left\{\sum_{\bd,n\ge 0}
A_{\bd,n}q^{\bd}\hbar^n\in \C[\hbar][\![q]\!]\;\Big|\; 
\begin{array}{l}
A_{0,n}=0 \text{ for } n>0, \\
|A_{\bd,n}|\le B C^{|\bd|} |\bd|^n (\exists B,C>0). 
\end{array}
\right\} 
\]
Let $\tFH$ be the following submodule of $\FH$:
\[
\tFH=\mO_{\rm small}^{\hbar}\langle P_1,\dots,P_r\rangle\Delta\subset\FH.
\] 
Then we have 
$\tFH\cong\mO_{\rm small}^{\hbar}\langle \hbar\partial \rangle 
/ \ov{I}_{\rm poly}$, 
where $\ov{I}_{\rm poly}$ is the closure of $I_{\rm poly}$ in 
$\mO_{\rm small}^{\hbar}\langle \hbar\partial \rangle$ with respect to 
the $q$-adic topology. 
Then we have a surjection 
\begin{equation}
\label{eq:surjfromMell}
M_{\rm Mell}\otimes_{\C[\hbar,q^{\pm}]} \mO^{\hbar}_{\rm small}[q^{-1}]
=\left(\mO^{\hbar}_{\rm small}\langle \hbar\partial\rangle
/\mO^\hbar_{\rm small} I_{\rm poly}\right)[q^{-1}] 
\longrightarrow \tFH[q^{-1}]
\end{equation}

\begin{proposition}
\label{prop:EFCdefinedovermO} 
Let 
$\{T_i(y_1,\dots ,y_r)\}_{i=0}^s$ be homogeneous polynomials  
such that 
$\{T_i(p_1,\dots,p_r)\}_{i=0}^s$ forms a basis of $H^*(X)$. 
Then, $\{T_i(P)\Delta\}_{i=1}^s$ forms a free basis of $\tFH$ 
as an $\mO^{\hbar}_{\rm small}$-module. 
Moreover, we have 
$\FH = \tFH\otimes_{\mO_{\rm small}^\hbar}\C[\hbar][\![q]\!]$. 
\end{proposition} 
\begin{proof}
In \cite{iritani-EFC}, section 4.4, we showed that 
$\{T_i(P)\Delta\}_{i=0}^s$ forms a free basis of 
$\FH$ as a $\C[\hbar][\![q]\!]$-module. 
Here, we will show this is also a basis of $\tFH$ 
as a $\mO_{\rm small}^\hbar$-module. 
It suffices to show that the connection matrix 
$\Omega_a=(\Omega_{a;ij})$ 
defined by 
\[
P_a T_j(P)\Delta =\sum_{i=0}^s \Omega_{a;ij}(q,\hbar)T_i(P)\Delta 
\]
has entries in $\mO_{\rm small}^\hbar$. 
We use a method similar to the proof in section 
\ref{subsect:powerserieswithest}. 
By \cite{iritani-EFC}, section 4.3,  
there exists a map 
$\Xi\colon \FH\rightarrow H^*(X)\otimes\C[\hbar,\hbar^{-1}][\![q]\!]$ 
and a pairing $(\cdot,\cdot)\colon 
FH^{S^1}_* \times \FH\rightarrow \C[\hbar][\![q]\!]$
such that 
\begin{equation}
\label{eq:EFCsolandpairing}
\Xi(P_a \alpha)=(\hbar\partial_a +p_a)\Xi(\alpha),\qquad
(\ov{\alpha},\beta)=\int_X \ov{\Xi(\alpha)}\cup \Xi(\beta), 
\end{equation}
for $\alpha,\beta\in \FH$. 
Here, the equivariant Floer homology $FH^{S^1}_*$ is the Poincar\'{e}
dual theory of $\FH$ with a $\C[\![q]\!]$-module isomorphism 
$\ov{\phantom{A}}\colon \FH\cong FH^{S^1}_*$ such that 
$\ov{\hbar\alpha}=-\hbar\ov{\alpha}$. 
The operator $\ov{\phantom{A}}$ also acts on 
$H^*(X)[\hbar,\hbar^{-1}][\![q]\!]$ by changing the sign of $\hbar$. 
The function $\tI(q,\hbar)=\Xi(\Delta)$ is written as  
\begin{equation}
\label{eq:I-function}
\tI(q,\hbar)
=\sum_{\bd\ge 0} q^{\bd}\prod_{i=1}^{r+N} 
\frac{\prod_{k=1+\pair{\td_i}{\bd}}^\infty (\td_i+k\hbar)} 
     {\prod_{k=1}^\infty (\td_i+k\hbar)}.
\end{equation}
Here, $I(q,\hbar)=e^{p\log q/\hbar}\tI(q,\hbar)$ is frequently 
referred to as $I$-function. 
We write $\tI(q,\hbar)=\sum_{\bd\ge 0}I_\bd(\hbar) q^{\bd}$. 
Take $C_1>0$ such that 
$\sup_{|\hbar|=1,k\ge 1,i}
(\|\td_i\|, \|1+\td_i/(k\hbar)\|, \|(1+\td_i/(k\hbar))^{-1}\| )
\le C_1$. 
Then we have for $|\hbar|=1$
\begin{align*}
\|I_{\bd}(\hbar)\| 
& \le 
\frac{\prod_{\pair{\td_i}{\bd}<0} |\pair{\td_i}{\bd}|!} 
     {\prod_{\pair{\td_i}{\bd}>0} \pair{\td_i}{\bd}!} 
C_1^{\sum_{i=1}^{r+N} |\pair{\td_i}{\bd}|}  \\
& \le C_2^{|\bd|}
\begin{cases}
1/\pair{c_1(X)}{\bd}! & \text{ if }\pair{c_1(X)}{\bd}\ge 0, \\
|\pair{c_1(X)}{\bd}|! & \text{ if }\pair{c_1(X)}{\bd}\le 0. 
\end{cases}
\end{align*}
for some $C_2>0$. 
Here, we used $c_1(X)=\sum_{i=1}^{r+N} \td_i$. 
Let $\mO_{\rm deg}^{\hbar,\hbar^{-1}}$ be the set of 
power series $\sum_{\bd} A_{\bd}(\hbar)q^{\bd}$ 
in $\C[\hbar,\hbar^{-1}][\![q]\!]$ satisfying 
\[
\sup_{|\hbar|=1}|A_{\bd}(\hbar)|\le B C^{\bd} 
\begin{cases}
1/\pair{c_1(X)}{\bd}! & \text{ if }\pair{c_1(X)}{\bd}\ge 0, \\
|\pair{c_1(X)}{\bd}|! & \text{ if }\pair{c_1(X)}{\bd}\le 0, 
\end{cases}
\] 
for some $B,C>0$. ($B$, $C$ depend on each element.)  
As in section \ref{subsect:powerserieswithest}, 
we can prove that $\mO_{\rm deg}^{\hbar,\hbar^{-1}}$ is a subring of 
$\C[\hbar,\hbar^{-1}][\![q]\!]$. 
Each component of $\tI(q,\hbar)$ belongs to the ring 
$\mO_{\rm deg}^{\hbar,\hbar^{-1}}$. 
Set $\eta_{kj}=(\ov{T_k(P)\Delta},T_j(P)\Delta)$.
By using (\ref{eq:EFCsolandpairing}), we have 
\begin{align*}
\eta_{kj}&=
\int_X T_k(-\hbar\partial+p)\tI(q,-\hbar) \cup 
       T_j(\hbar\partial+p)\tI(q,\hbar) \\
\sum_{i=0}^s \eta_{ki}\Omega_{a;ij}
&=(\ov{T_k(P)\Delta}, P_aT_j(P)\Delta)  \\
&=\int_X T_k(-\hbar\partial+p)\tI(q,-\hbar) \cup 
       (\hbar\partial_a+p_a) T_j(\hbar\partial+p)\tI(q,\hbar).
\end{align*}
Therefore, $\eta_{kj}$ and $\sum_{i}\eta_{ki}\Omega_{a;ij}$ 
are contained in $\mO_{\rm deg}^{\hbar,\hbar^{-1}}\cap \C[\hbar][\![q]\!]$. 
Moreover, these elements are homogeneous 
with respect to the grading $\deg q^{\bd}=2\pair{c_1(X)}{\bd}$ 
and $\deg \hbar=2$. 
The homogeneity means that the coefficient of $q^{\bd}\hbar^n$ 
is non-vanishing only when $\pair{c_1(X)}{\bd}+n$ equals a given constant. 
Therefore, by Cauchy's residue theorem, we can see that 
$\eta_{kj}$ and $\sum_{i}\eta_{ki}\Omega_{a;ij}$ are also contained in 
$\mO_{\rm small}^{\hbar}$. 
Because $\lim_{q\to 0}\eta_{kj}=\int_X T_k(p)\cup T_j(p)$ is an invertible 
matrix and $\mO_{\rm small}^{\hbar}$ is a local ring, 
we have $\Omega_{a;ij}\in \mO_{\rm small}^{\hbar}$.
\end{proof}

\subsection{Generalized mirror transformation} 
We can describe how to reconstruct the quantum $D$-module $QDM^*(X)$ 
from $\FH$ as follows: 
\begin{theorem}[{\cite{jinzenji2}}, {\cite{iritani-EFC}} Theorem 4.9, 5.4]
\label{thm:generalizedmirrortransformation}
There exists a formal embedding 
$\emb\colon  (\C^{r},0)\rightarrow (\C^{s+1},0)$ 
and an isomorphism of $D$-modules 
$\Phi_{\emb}\colon \emb^*(QDM(X))\cong \FH$. 
The map $\emb$ is given by equations of the form 
\[
\hatt_0=F_0(q), \hq_1=q_1\exp(F_1(q)),\dots,\hq_r=q_r\exp(F_r(q)),
\hatt_{r+1}=F_{r+1}(q),\dots ,\hatt_{s}=F_s(q) 
\]
for some $F_i(q)\in \C[\![q]\!]$, $F_i(0)=0$ 
and ${\Phi_{\emb}}|_{q=0}$ is determined by the canonical isomorphism 
$H^*(X)\otimes\C[\hbar]\cong \FH/\sum_{a=1}^r q_a \FH$. 
Moreover, we can reconstruct $QDM^*(X)$ from $\FH$ by the following 
steps: 

{\rm (i)} Take a free basis $\{T_0,\dots,T_s\}$ 
of $\FH$ as a $\C[\hbar][\![q]\!]$-module as in 
Proposition \ref{prop:EFCdefinedovermO}
and calculate a connection matrix $\Omega_a$ 
defined by $P_aT_j=\sum_{i=0}^s \Omega_{a;ij}T_i$. 

{\rm (ii)} Find a gauge transformation $g$ such that $g|_{q=0}=\id$ 
and that the new connection matrix $\hOmega_a$ is $\hbar$-independent, 
where $\hOmega_a = g^{-1} \Omega_a g + g^{-1} \hbar \partial_a g$.  

{\rm (iii)} Solve for matrix-valued functions 
$\hOmega_j(t_0,q,t_{r+1},\dots,t_s)$ 
$(0\le j\le s)$ from $\hOmega_a(q)$ $(1\le a\le r)$, 
where $\hOmega_0=\id$.  
This procedure will be reviewed 
in Proposition \ref{prop:reconstructionandconvergence}.

{\rm (iv)} Find a new coordinate system 
$(\hatt_0,\hq_1,\dots,\hq_r,\hatt_{r+1},\dots,\hatt_s)$ 
of the form $\hatt_0=t_0+F_0(q)$, 
$\hq_a=q_a\exp(F_a(q))$ $(1\le a\le r)$, 
$\hatt_j=t_j+F_j(q)$ $(r+1\le j \le s)$
such that $\hOmega_{\hat{j}}=\sum_{i=0}^s (\partial t_i/\partial \hatt_j) 
\hOmega_i$ satisfies $\hOmega_{\hat{j}}(1)=p_j$.   
\end{theorem} 

We study these four steps from an analytic point of view.  
By Proposition \ref{prop:EFCdefinedovermO}, 
the connection matrix $\Omega_a(q)$ in step (i) has 
its matrix elements in $\mO^\hbar_{\rm small}$. 
Then by Proposition \ref{prop:gaugefixing}, the connection matrix 
$\hOmega_a$ in step (ii) becomes a convergent function of $q$. 
In Proposition \ref{prop:reconstructionandconvergence}, 
we will see that the reconstruction step 
(iii) preserves the convergence. 
The last step (iv) can be done in the convergent category, 
therefore $QH^*(X)$ has convergent structure constants. 
Summarizing, 
\begin{theorem}
\label{thm:analyticityofgenmirtrans}
The quantum $D$-module $QDM^*(X)$ of a toric variety 
is defined over convergent power series.  
The embedding $\emb\colon (\C^r,0)\rightarrow (\C^{s+1},0)$ in 
Theorem \ref{thm:generalizedmirrortransformation} is complex analytic 
and 
$\tFH\cong \emb^*(QDM^*_{\rm an}(X))
\otimes_{\mO_{\rm an}^\hbar} \mO^{\hbar}_{\rm small}$, 
where $QDM^*_{\rm an}(X)=(H^*(X)\otimes\mO^\hbar_{\rm an},\nabla^{\hbar})$ 
and $\mO^\hbar_{\rm an}=
\{f(\hbar,x)\in \C\{\hbar,x\}\;|\; f|_{q=0}\text{ is constant.}\}$. 
\end{theorem}

In order to state the compatibility of reconstruction with convergence, 
we consider the following general situation. 
Consider a flat connection 
$\nabla^\hbar=\hbar d+\sum_{a=1}^r \Omega_a(q) d q_a/q_a$ 
of the bundle $H^*(X)\times U\rightarrow U$ 
regular singular along $q_1\cdots q_r=0$. 
Here, $U$ is a neighborhood of $0$ in $\C^r$ and  
$\Omega_a(q)$ is an $\hbar$-independent holomorphic function on $U$. 
\begin{proposition}
\label{prop:reconstructionandconvergence}
Assume that $H^*(X)$ is generated by $1$ 
under the action of residue matrices $p_a:=\Omega_a(0)$. 
For a set of vectors $v_1,\dots,v_l$ in $H^*(X)$, 
$\nabla^{\hbar}$ can be uniquely extended to a flat connection 
$\widetilde{\nabla}^{\hbar}=
\hbar d+\sum_{a=1}^r \Omega_a(q,t) d q_a/q_a 
+\sum_{j=1}^l \Phi_j(q,t) dt_j$ 
defined on a neighborhood $U'$ of $0$ in $\C^{r+l}$ 
under the condition that $\Phi_j(q,t)1=v_j$.  
\end{proposition}
\begin{proof}
In \cite{iritani-genmir}, Theorem 4.9, we proved that 
$\widetilde{\nabla}^\hbar$ can be reconstructed uniquely 
as a formal connection. 
Here, we prove that this is convergent. 
Because we can extend $\nabla^{\hbar}$ to the $t_j$-direction 
independently for different $t_j$'s, 
we can assume that $l=1$. 
Flatness of the connection $\widetilde{\nabla}^{\hbar}$ implies 
$\partial_t \Omega_a = q_a\partial_{q_a}\Phi$ and 
$[\Omega_a,\Phi]=0$. 
By taking a smaller $U$ if necessary, 
we can assume that $H^*(X)$ is generated by $1$ under the 
action of $\Omega_a(q)$ for each $q\in U$.  
Then, we have a surjection 
$\C[\Omega_1(q),\dots,\Omega_r(q)]\rightarrow H^*(X)$ 
and this determines a product $*_q$ on $H^*(X)$ for each $q\in U$. 
On the locus $\{t=0\}$, $\Phi(q,0)$ must be a multiplication  
$(v_1*_q)$ by $v_1$ because it commutes with $\Omega_a(q)$.  
Therefore, $\Phi(q,0)$ becomes holomorphic on $U$. 
On the locus $\{q=0\}$, it is easy to see that 
$\Omega_a(0,t)=p_a$ and $\Phi(0,t)=v_1\cup$. 
Expand $\Phi$ and $\Omega_a$ as 
$\Phi=\sum_{\bd,m\ge 0} \Phi_{\bd,m}q^\bd t^m$ and 
$\Omega_a=\sum_{\bd,m\ge 0} \Omega_{a;\bd,m} q^\bd t^m$. 
Note that $\Phi_{0,m}=\delta_{0m}(v_1\cup)$ and 
$\Omega_{a;0,m}=\delta_{0m}p_a$. 
Then by $\partial_t \Omega_a = q_a\partial_{q_a}\Phi$, we have 
\begin{equation}
\label{eq:phi2omega}
\Omega_{a;\bd,m+1}=\frac{d_a}{m+1} \Phi_{\bd,m}.
\end{equation}
Let $\{T_0,\dots,T_s\}$ be a basis of $H^*(X)$. 
By the assumption, we can write 
\[
T_i=\sum_{l,a_1,\dots, a_l} A_{a_1,\dots,a_l}^{(i)} p_{a_1}\cdots p_{a_l}(1)
\]
for some $A_{a_1,\dots,a_l}^{i}\in \C$. 
Then by $\Phi_{\bd,m+1}(1)=0$, we have 
\begin{align}
\label{eq:phim+1}
\Phi_{\bd,m+1}T_i &=
\sum_{l,a_1,\dots, a_l} A_{a_1,\dots,a_l}^{(i)} 
\sum_{k=1}^l p_{a_1}\cdots [\Phi_{\bd,m+1},p_{a_k}]\cdots p_{a_l}(1) \\
\nonumber 
&= \sum_{l,a_1,\dots, a_l} A_{a_1,\dots,a_l}^{(i)} 
\sum_{k=1}^l p_{a_1}\cdots 
\left(\sum_{\bd_1>0}\sum_{j=0}^{m+1}
[\Omega_{a_k;\bd_1,j},\Phi_{\bd-\bd_1,m+1-j}]\right)\cdots p_{a_l}(1) 
\end{align}
where we used $[\Phi,\Omega_{a_k}]=0$. 
Assume by induction that we know $\Omega_{a;\bd',m}$ and 
$\Phi_{\bd',m}$ for all $|\bd'|< \ovd$ 
and all $m\ge 0$. 
Set $\omega_{d,m}=\max_a \sum_{|\bd|=d} \|\Omega_{a;\bd,m}\|$ 
and $\phi_{d,m}=\sum_{|\bd|=d} \|\Phi_{\bd,m}\|$. 
Assume also that there exist constants $A,C,M>1$ such that 
\begin{equation}
\label{eq:estimateofomegaandphi}
\omega_{d,m}\le A\frac{C^{d+m}}{(d+1)^M}\frac{d^m}{m!}, \quad 
\phi_{d,m} \le A\frac{C^{d+m}}{(d+1)^M}\frac{d^m}{m!}
\end{equation} 
holds for all $d<\ovd$ and all $m\ge 0$. 
We must choose $A,C,M$ so that this estimate is valid for $m=0$ and 
all $d\ge 0$.  
Take $\bd$ such that $|\bd|=\ovd$. 
First we can solve for $\Omega_{a;\bd,1}$ by (\ref{eq:phi2omega}). 
Then we solve for $\Phi_{\bd,1}$ by (\ref{eq:phim+1}). 
Next we solve for $\Omega_{a;\bd,2}$ by (\ref{eq:phi2omega}) 
and repeat this process. 
Assume that the estimate (\ref{eq:estimateofomegaandphi}) 
holds for $d=\ovd$ and up to $m$ ($m\ge 0$). 
We have by (\ref{eq:phi2omega}), 
\[
\omega_{\ovd,m+1}\le \frac{\ovd}{m+1} \phi_{\ovd,m} 
\le A \frac{C^{\ovd+m}}{(\ovd+1)^M} \frac{\ovd^{m+1}}{(m+1)!}.
\]
Next we have by (\ref{eq:phim+1}), 
\begin{align}
\label{eq:estofphim+1}
\phi_{\ovd,m+1}
&\le B_1\sum_{0<i<\ovd}\sum_{j=0}^{m+1}\omega_{i,j}\phi_{\ovd-i,m+1-j} 
+B_1\omega_{\ovd,m+1} \|(v_1\cup)\| \\
\nonumber
&\le B_1 A^2 \frac{C^{\ovd+m+1}}{(\ovd+1)^M} \frac{\ovd^{m+1}}{(m+1)!}
\sum_{i=1}^{\ovd-1} 
\frac{(\ovd+1)^M}{(i+1)^M(\ovd-i+1)^M} 
+B_2A\frac{C^{\ovd+m}}{(\ovd+1)^M} \frac{\ovd^{m+1}}{(m+1)!} \\
\nonumber
&\le A \frac{C^{\ovd+m+1}}{(\ovd+1)^M} \frac{\ovd^{m+1}}{(m+1)!}
\{ B_1A\varepsilon_1(M)+B_2/C \}
\end{align}
where $B_1,B_2 >0$ are constants determined only by $(H^*(X),\cup)$ 
and $\varepsilon_1(M)$ is a function defined in the proof of 
Proposition \ref{prop:gaugefixing}. 
To complete the induction, we need to choose $C$ and $M$ 
carefully. 
Because $\lim_{M\to \infty}\varepsilon_1(M)=0$, for 
sufficiently large $M$, we have $B_1A\varepsilon_1(M)\le 1/2$. 
Next we choose $C$ sufficiently large so that 
the estimate (\ref{eq:estimateofomegaandphi}) is valid 
for $m=0$ and all $d\ge 0$ and that $C\ge 2 B_2$. 
Then by (\ref{eq:estofphim+1}), we have the desired estimate for 
$\phi_{\ovd,m+1}$. 
After we obtain the estimate (\ref{eq:estimateofomegaandphi})
for all $d,m\ge 0$, 
we can see that $\Omega_a(q,t)$ and $\Phi(q,t)$ are convergent  
because 
\[
\sum_{d\ge 0} \sum_{m\ge 0} C^{d+m} \frac{d^m}{m!} x^d y^m 
=\frac{1}{1-Cx\exp(Cy)}
\]
is holomorphic around $(x,y)=(0,0)$. 
\end{proof}
When we apply this proposition to the dual Givental connections, we have
\begin{corollary}
\label{cor:convergenceofbig}
Let $X$ be a smooth projective variety. 
If $H^*(X)$ is generated by $H^2(X)$ and 
if the small quantum cohomology of $X$ has convergent structure constants, 
so does the big quantum cohomology $QH^*(X)$.  
\end{corollary}

\subsection{Characteristic variety and semisimplicity}
We study the characteristic variety of the Mellin system 
and proves the semisimplicity. 
For a differential operator $f(q,\hbar\partial,\hbar)$, 
we define its principal symbol as $\sigma(f)=f(q,\cp,0)$.  
Here, $\cp_1,\dots,\cp_r$ are conjugate variables.  
We also define 
\begin{align*}
\sigma(M_{\rm Mell})&=\C[q^{\pm},\cp]/\sigma(I_{\rm Mell}),\\
\sigma(\tFH)&=\C\{q\}[\cp]\otimes_{\sigma} \tFH
\cong \C\{q\}[\cp]/\sigma(\ov{I}_{\rm poly}). 
\end{align*}
The characteristic varieties of $\Ch(M_{\rm Mell})$ and 
$\Ch(\tFH)$ are defined as analytic spectra of $\sigma(M_{\rm Mell})$ 
and $\sigma(\tFH)$ respectively. 
They are (germs of) analytic subvarieties of $(\C^*)^r\times\C^r=T^*(\C^*)^r$. 
By the surjection (\ref{eq:surjfromMell}), we have an embedding 
\[
\Ch(\tFH)|_{q_1\cdots q_r\neq 0} \longrightarrow 
\Ch(M_{\rm Mell})
\]
The characteristic variety of $M_{\rm Mell}$ is 
the zero set of $\sigma(\mP_{\bd})$ for all $\bd$.  
Because there exists $\bd$ such that $\pair{\td_i}{\bd}>0$ for all $i$, 
we can see from $\sigma(\mP_{\bd})=0$ that 
$\sum_{b=1}^r m_{ib}\cp_b\neq 0$ if $q\in (\C^*)^r$. 
(Here, we need the compactness of $X$.) 
Thus, 
\[
\Ch(M_{\rm Mell})=\Bigl\{(q,\cp)\in (\C^*)^r\times \C^r \;\Big|\; 
\sum_{b=1}^r m_{ib}\cp_b\neq 0, \ 
q_a=\prod_{i=1}^{r+N} \Bigl(\sum_{b=1}^r m_{ib}\cp_b\Bigr)^{m_{ia}}
 \Bigr\}.
\]
Let $\Crit(F_q)$ be the set of critical points of $F_q$ in $Y_q$. 
This forms a family over $(\C^*)^r$. 
For $x\in \Crit(F_q)$, the cotangent vector $d_xF$ in $T^*_xY$ equals  
$\pi^*(\sum_{a=1}^r \cp_a d\log q_a)$ for some $\cp_a$.  
Hence we have a map $dF\colon \bigcup_{q\in (\C^*)^r} \Crit(F_q)
\rightarrow T^*(\C^*)^r$, 
$x\mapsto \sum_{a=1}^r \cp_a d\log q_a$. 

\begin{lemma}
\label{lem:criticalfamilyischaracteristic}
$dF\colon \bigcup_{q\in (\C^*)^r} \Crit(F_q)
\cong \Ch(M_{\rm Mell}) \subset T^*(\C^*)^r$.
\end{lemma}
\begin{proof}
By $dF=\pi^*(\sum_{a=1}^r \cp_a d\log q_a)$ on $\bigcup_q \Crit(F_q)$, 
we see that $\mc_i=\sum_{a=1}^r m_{ia}\cp_a$ and so  
$(q,\cp)\in \Ch(M_{\rm Mell})$. 
The inverse is given by $\mc_i=\sum_{a=1}^r m_{ia}\cp_a$.  
\end{proof}

\begin{proposition}
\label{prop:coveringdegree}
{\rm (i)} The projection $\Ch(M_{\rm Mell})\rightarrow (\C^*)^r$ 
is a ramified covering of degree $N!\Vol(\Conv(\Sigma))$, 
where $\Conv(\Sigma)$ is the convex hull of all the primitive generators 
$\vx_1,\dots, \vx_{r+N}$ 
of 1-dimensional cones in the fan $\Sigma$ defining $X$.

{\rm (ii)} There exist exactly $\dim_{\C}H^*(X)$ branches of the 
covering $\Ch(M_{\rm Mell})\rightarrow (\C^*)^r$ corresponding to 
the subcovering $\Ch(\tFH)|_{q_1\cdots q_r\neq 0} 
\rightarrow (\C^*)^r$. 
These branches are characterized by the condition that 
$\cp_a \to 0$ as $q\to 0$.  
\end{proposition} 

\begin{proof}
By Kushnirenko's theorem \cite{kushnirenko}, 
the dimension of Jacobi ring 
$\C[s_1,\dots,s_N]/\langle
\partial_{s} F_q \rangle$ 
equals $N!$ times the volume of the Newton polytope of 
$F_q=\sum_{i=1}^{r+N} (\prod_{a} q_a^{l_{ia}}) s^{\vx_i}$. 
On the other hand, the projection 
$\pi\colon \Ch(M_{\rm Mell})\rightarrow (\C^*)^r$ 
is a submersion at generic $\cp$. 
(In fact, $d\pi=
\partial \log q_a/\partial \cp_b=\sum_{i}m_{ia}\mc_i^{-1}m_{ib}$ 
is positive definite when $\mc_i>0$.)  
Hence we obtain (i) by  Lemma \ref{lem:criticalfamilyischaracteristic}.
Each branch of $\Ch(\tFH)$ corresponds to 
a branch of simultaneous eigenvalues 
$(\cp_1,\dots,\cp_r)$ of the connection matrices 
$(\Omega_1|_{\hbar=0},\dots,\Omega_r|_{\hbar=0})$. 
Because $\Omega_a|_{q=0}$ is a nilpotent matrix (cup product by $p_a$), 
we see that $\cp_a\to 0$ as $q\to 0$. 
On the other hand, we have $\Ch(M_{\rm Mell})\cong 
\Ch(FH_0)|_{q_1\cdots q_r\neq 0}$ 
and the zero-fiber of $\Ch(FH_0)$ is the spectrum of 
$
\sigma(FH_0)/\sum_{a=1}^r q_a \sigma(FH_0)\cong 
\sigma(FH_0/\frm FH_0)\cong \sigma(\widehat{FH}_0/\frm \widehat{FH}_0) 
\cong H^*(X)
$. 
Therefore, there exist only $\dim_{\C}H^*(X)$ branches 
converging to $\cp=0$ at $q=0$. 
\end{proof}

By Theorem \ref{thm:analyticityofgenmirtrans}, we see that 
$\emb^*(\Ch(QDM_{\rm an}^*(X)))\cong \Ch(\tFH)$. 
Therefore, the generic fiber of the characteristic variety of 
$QDM^*_{\rm an}(X)$  
consists of  $\dim_{\C} H^*(X)$ distinct points. 
This means that 
\begin{corollary}
\label{cor:semisimplicity}
The quantum cohomology of a smooth, projective toric variety is 
generically semisimple. 
(See \cite{manin}, Part I, section 3 for the definition of semisimplicity.)
\end{corollary} 
As remarked in the introduction, this corollary 
together with Kawamata's result \cite{kawamata} 
shows that Bayer and Manin's modified Dubrovin's conjecture 
\cite{bayer-manin} holds for toric varieties. 

The proposition below will clarify the role of the nef condition 
for $c_1(X)$.   
\begin{proposition}
\label{prop:nefcharacterization}
The following conditions are equivalent. 

{\rm (i)} $c_1(X)$ is nef. 

{\rm (ii)} $N!\Vol(\Conv(\Sigma))=\dim_{\C}H^*(X)$.

{\rm (iii)} $\Ch(\tFH)|_{q_1\cdots q_r\neq 0}\cong
\Ch(M_{\rm Mell})$ on a small neighborhood of $q=0$.  

{\rm (iv)} $\tFH\cong FH_0\otimes_{\C[\hbar,q]}\mO^\hbar_{\rm small}
(\cong 
\mO^{\hbar}_{\rm small} \langle \hbar\partial \rangle /
\mO^{\hbar}_{\rm small} I_{\rm poly})$. 

{\rm (v)} $\FH\cong FH_0\otimes_{\C[\hbar,q]}\C[\hbar][\![q]\!]
(\cong 
\C[\hbar][\![q]\!] \langle \hbar\partial \rangle /
\C[\hbar][\![q]\!] I_{\rm poly})$. 

{\rm (vi)} The generalized mirror transformation can be done 
using only convergent power series i.e.
$\Omega_a(q,\hbar)$ and $g(q,\hbar)$ in the step {\rm (i)} and {\rm (ii)} of 
Theorem \ref{thm:generalizedmirrortransformation} 
are convergent functions of $q$ and $\hbar$. 
\end{proposition}
\begin{proof}
(i) $\Leftrightarrow$ (ii): 
The condition $c_1(X)\ge 0$ is equivalent to 
that every primitive generator $\vx_i$ of one-dimensional cones 
lies in the boundary of $\Conv(\vx_1,\dots,\vx_{r+N})$. 
The number of top dimensional cones in $\Sigma$ equals 
$\dim_{\C}H^*(X)$ and each top dimensional cone has volume $1/N!$. 
Hence we obtain the equivalence. 

(ii) $\Leftrightarrow$ (iii): 
This follows from Proposition \ref{prop:coveringdegree}.

(iii) $\Rightarrow$ (iv): 
Set $M:=\mO^{\hbar}_{\rm small}\langle \hbar\partial \rangle /
\mO^{\hbar}_{\rm small} I_{\rm poly}$. 
Let $[f(q,\hbar\partial,\hbar)]$ be an element of 
the kernel of the natural surjection $M\rightarrow \tFH$. 
We can assume that $f(q,\hbar\partial,\hbar)$ is homogeneous. 
By the assumption, we can see that 
there exists an integer $k>0$ such that 
$(q_1\cdots q_r)^k f(q,\cp,0)$ is contained in 
$\sigma(\mO^{\hbar}_{\rm small}I_{\rm poly})$. 
Because $\sigma(\mO^{\hbar}_{\rm small}I_{\rm poly})=
\sum_{\bd} \C\{q\}[q^{-1}]\sigma(\mP_{\bd}) \cap \C\{q\}$, 
we can see that $f(q,\cp,0)$ is also in 
$\sigma(\mO^{\hbar}_{\rm poly} I_{\rm poly})$. 
Therefore, there exists some $f_1(q,\hbar\partial,\hbar)$ 
in $\mO^\hbar_{\rm small}\langle \hbar\partial \rangle$
such that $f(q,\hbar\partial,\hbar)\equiv \hbar f_1(q,\hbar\partial,\hbar)
\mod \mO^\hbar_{\rm small}I_{\rm poly}$ and $\deg f_1=\deg f-2$. 
Then, $[\hbar f_1]$ is in the kernel of $M\rightarrow \tFH$. 
Because the multiplication by $\hbar$  is injective in $\tFH$, 
$[f_1]$ also lies in the kernel. 
Repeating this, we have 
$f\equiv \hbar^n f_n \mod \mO^\hbar_{\rm small}I_{\rm poly}$ 
and $\deg f_n<0$ for sufficiently large $n$. 
Because we have already seen that (i) is equivalent to (iii), 
$c_1(X)$ is nef and so $\deg q_a\ge 0$. 
Thus $f_n$ must be zero and $[f]=0$ in $M$. 

(iv) $\Rightarrow$ (v): Take the tensor product with $\C[\hbar][\![q]\!]$. 

(v) $\Rightarrow$ (ii): 
Because $\FH$ is a free $\C[\hbar][\![q]\!]$-module of rank 
$\dim_{\C} H^*(X)$, 
we see that $\sigma(FH_0)\otimes_{\C[q]} \C[\![q]\!]$ is a
free $\C[\![q]\!]$-module of rank $\dim_{\C} H^*(X)$. 
On the other hand, by Proposition \ref{prop:coveringdegree}  
$\sigma(FH_0)[q^{-1}]\cong \sigma(M_{\rm Mell})$  
becomes a free module of rank $N!\Vol(\Conv(\Sigma))$ 
when restricted to the complement of the ramification locus. 
Therefore, we must have $\dim_{\C} H^*(X)=N!\Vol(\Conv(\Sigma))$. 

(i) $\Leftrightarrow$ (v): 
We can see from the proof of Proposition \ref{prop:EFCdefinedovermO}
that the $I$-function in (\ref{eq:I-function}) 
is convergent on $|\hbar|=1$ if and only if $c_1(X)$ is nef. 
The connection matrix $\Omega_a$ is determined from $I$-function by 
(\ref{eq:EFCsolandpairing}). 
The gauge transformation $g$ 
can be found by applying Birkhoff factorization of the loop group 
$LGL(s+1,\C)$ 
to the loop $\hbar\mapsto (T_i(\hbar\partial)I(q,\hbar))$ 
(see \cite{guest,iritani-EFC,iritani-genmir}). 
\end{proof}

\subsection{Example}
We study the case of Hirzeburch surface $\F_n=\Proj(\mO_{\Proj^1}
\oplus \mO_{\Proj^1}(n))$. 
This is given by the fan $\Sigma$ whose one-dimensional cones are 
generated by $\vx_1,\dots,\vx_4$:    
\[
\vx_1=(1,0),\quad \vx_2=(0,1),\quad \vx_3=(-1,n),\quad 
\vx_4=(0,-1). 
\]
We have $\dim_{\C}H^*(\F_n)=4$ and 
$2!\Vol(\Conv(\Sigma))=\max(4,n+2)$.
Let $\td_1,\dots,\td_4$ be the classes of corresponding toric divisors. 
Then, $p_1:=\td_1$ and $p_2:=\td_4$ form a nef integral basis 
and we have $\td_2=p_2-np_1$ and $\td_3=p_1$. 
The relation of classical cohomology ring is given by 
$p_1^2=0$ and $p_2^2=np_1p_2$. 
The corresponding Mellin system is generated by the following 
two differential operators:  
\[
\mP_{(1,0)}=q_1 \tdD_2(\tdD_2-\hbar)\cdots(\tdD_2-(n-1)\hbar)- \tdD_1\tdD_3,
\quad 
\mP_{(0,1)}=q_2-\tdD_4\tdD_2,
\]
where $\tdD_1=\tdD_3=\hbar\partial_1$, 
$\tdD_2=\hbar\partial_2-n\hbar\partial_1$, 
$\tdD_4=\hbar\partial_2$. 
\begin{gather*}
FH_0\cong \C\langle q_1,q_2,\hbar\partial_1,\hbar\partial_2,\hbar\rangle /
\langle \mP_{(1,0)},\mP_{(0,1)} \rangle,\quad  
M_{\rm Mell} \cong FH_0[q_1^{-1},q_2^{-2}]\\
\FH\cong \C[\hbar][\![q_1,q_2]\!]\langle\hbar\partial_1,\hbar\partial_2\rangle 
\Big/
\ov{\langle \mP_{(1,0)},\mP_{(0,1)} \rangle}. 
\end{gather*}
\subsubsection{Fano case $(n=0,1)$}
In this case, $FH_0$ is freely generated by $[1], [\hbar\partial_1], 
[\hbar\partial_2], [(\hbar\partial_2)^2]$ over 
$\C[\hbar,q]$. 
The small quantum $D$-module of $\F_n$ is defined over 
the polynomial ring $\C[\hbar,q]$ 
and is isomorphic to $FH_0$. 
\subsubsection{Nef but  non-Fano case $(n=2)$}
In this case, $FH_0$ is not finitely generated over $\C[\hbar,q]$. 
However,   $FH_0[(1-4q_1)^{-1}]$ is freely generated by 
$[1], [\hbar\partial_1], 
[\hbar\partial_2], [(\hbar\partial_2)^2]$ over 
$\C[\hbar,q_1,q_2,(1-4q_1)^{-1}]$. 
For example, we can write $[(\hbar\partial_1)^2]$ as 
\[
[(\hbar\partial_1)^2]=\frac{2q_1q_2}{1-4q_1}
+\frac{2\hbar q_1}{1-4q_1}[\hbar\partial_1]
-\frac{\hbar q_1}{1-4 q_1}[\hbar\partial_2]
-\frac{q_1}{1-4 q_1}[(\hbar\partial_2)^2]
\]
The small quantum $D$-module is isomorphic to $FH_0$ 
at least on the region $|4 q_1|<1$ 
after a suitable coordinate change (see \cite{guest}). 
\subsubsection{Non-nef case $(n\ge 3)$}
In this case, $FH_0$ is not finitely generated over $\C[\hbar,q]$.  
Furthermore, $FH_0\otimes_{\C[\hbar,q]}\C[\hbar][\![q]\!]
\cong \C[\hbar][\![q]\!]\langle\hbar\partial \rangle/ 
\langle \mP_{(1,0)},\mP_{(0,1)} \rangle$ 
is not finitely generated over $\C[\hbar][\![q]\!]$ either. 
If it were finitely generated, then 
$\sigma(FH_0)\otimes_{\C[q]}\C[\![q]\!]$ would be also 
finitely generated over $\C[\![q]\!]$ and 
generated by $1,\cp_1,\cp_2,\cp_2^2$ by Nakayama's lemma. 
It contradicts that 
$\sigma(FH_0)[q_1^{-1},q_2^{-1}]\cong \sigma(M_{\rm Mell})$ 
is a free module of rank $n+2>4$. 
We must consider the $q$-adic closure of the left ideal 
$\langle\mP_{(1,0)}, \mP_{(0,1)}\rangle$ which is strictly bigger: 
\[
\langle \mP_{(1,0)},\mP_{(0,1)} \rangle \subsetneqq
\ov{\langle \mP_{(1,0)},\mP_{(0,1)} \rangle }\quad 
\text{ in }  \C[\hbar][\![q]\!]\langle \hbar\partial \rangle.  
\] 
Our $\FH$ is a free $\C[\hbar][\![q]\!]$-module of rank 4 
and isomorphic to the restriction of the big quantum $D$-module 
to some non-linear subvariety of $H^*(\F_n)$.   
On the other hand, $M_{\rm Mell}$ is of rank $n+2$ over 
$\C[\hbar,q^{\pm}]$.  
(In fact, it is generated by 
$[1],[\hbar\partial_1],[(\hbar\partial_1)^2], 
[\hbar\partial_2-n\hbar\partial_1],\dots, 
[(\hbar\partial_2-n\hbar\partial_1)^{n-1}]$.)  
The characteristic variety of $M_{\rm Mell}$ is given by 
\[
q_1(\cp_2-n\cp_1)^n=\cp_1^2,\quad q_2=(\cp_2-n\cp_1)\cp_2
\]
or equivalently,  
\[
n^2 q_1q_2^n=\cp_2^{n-2}(\cp_2^2-q_2)^2, \quad 
q_2=(\cp_2-n\cp_1)\cp_2.
\]
It has $n+2$ solutions $(\cp_1,\cp_2)$ (counted with multiplicity) 
for a given $(q_1,q_2)\in (\C^*)^2$. 
Out of $n+2$ solutions, we have four branches of solutions 
having the asymptotics: 
\[
\cp_1\sim (q_1(\sqrt{q_2})^n)^{1/2},\quad 
\cp_2\sim \sqrt{q_2}+\frac{n}{2} (q_1(\sqrt{q_2})^n)^{1/2}.
\]
These branches satisfy $(\cp_1,\cp_2)\to 0$ as $(q_1,q_2)\to 0$ 
and correspond to the characteristic variety of $\tFH$. 
The other $n-2$ branches have the asymptotics
\[\cp_1\sim -n^{\frac{n-4}{n-2}}q_1^{-\frac{1}{n-2}},\quad  
\cp_2\sim n^{\frac{2}{n-2}}q_1^{\frac{1}{n-2}} q_2
\] 
and diverge as $q_1\to 0$.

By changing the large radius limit $q\to 0$, 
we can construct a $D$-module of rank $n+2$ which is  
regular singular along $\hq_1\hq_2=0$. 
Set $\hq_1=q_1^{-1}$ and $\hq_2=q_2q_1$. 
In terms of $\hq_1,\hq_2$, differential operators of the Mellin system 
are written as 
\begin{gather*}
\mR_{(1,0)}=\hq_1\tdD_1\tdD_3-\tdD_2(\tdD_2-\hbar)
\cdots(\tdD_2-(n-1)\hbar),\quad
\mR_{(1,1)}=\hq_1\hq_2-\tdD_2\tdD_4, \\
\mR_{(0,1)}=\hq_2\tdD_2(\tdD_2-\hbar)\cdots (\tdD_2-(n-2)\hbar)-
\tdD_1\tdD_3\tdD_4.
\end{gather*}
Then, the $D$-module 
\[FH_0'=\C\langle\hq_1,\hq_2, \hbar\partial_1,\hbar\partial_2,\hbar\rangle /
\langle \mR_{(1,0)},\mR_{(1,1)},\mR_{(0,1)}\rangle
\]
is freely generated by $[1],[\hbar\partial_1],[(\hbar\partial_1)^2],
[\hbar\partial_2-n\hbar\partial_1], \cdots,
[(\hbar\partial_2-n\hbar\partial_1)^{n-1}]$ 
as a $\C[\hbar,\hq_1,\hq_2]$-module. 
This becomes an abstract quantum $D$-module in the sense of 
\cite{iritani-EFC,iritani-genmir}. 
Therefore, we can find a canonical frame and flat coordinates  
by reconstructing an $n+2$-dimensional base space 
(see \cite{iritani-genmir}).  
It would be interesting to study a geometric meaning of
this $D$-module.

\section{$R$-conjecture (Virasoro constraints)}
The Virasoro constraints are 
infinite dimensional symmetries of the (all genus, descendant) 
Gromov-Witten potential 
conjectured by Eguchi-Hori-Xiong \cite{eguchi-hori-xiong}. 
In \cite{givental-quadratic}, Givental showed that when 
the target $X$ has a torus action with only isolated fixed points 
and isolated one-dimensional orbits, 
the Virasoro conjecture is reduced to the {\it $R$-conjecture}. 
The $R$-conjecture is a conjecture about equivariant quantum cohomology 
with semisimple non-equivariant counterpart.  
It states that the asymptotic solution $R$ defined 
in the equivariant theory has non-equivariant limit. 
For a class of Fano toric varieties, the $R$-conjecture 
was proved by Givental \cite{givental-quadratic}. 
The $R$-conjecture has been proved for 
complete flag \cite{joe-kim}, Grassmannian \cite{bertram-ciocan-kim-twoproofs} 
and partial flag varieties \cite{bertram-ciocan-kim-nonabelian}. 

In this section, we explain the $R$-conjecture and 
prove it for any smooth projective toric variety. 
The proof of $R$-conjecture will be a first step 
for the understanding of the mirror oscillatory integral. 

\subsection{Equivariant quantum cohomology}
\label{sect:equivariantquantumcohomology}
Let $\T$ be an $l$-dimensional torus ($\cong (S^1)^l$) 
and $X$ be a $\T$-manifold. 
In a manner analogous to section \ref{sect:quantumD-modules}, 
we can define the equivariant quantum cohomology algebra 
$(QH^*_{\T}(X),*)$
which is a deformation of the ring structure of 
the ordinary equivariant cohomology $H^*_{\T}(X)$ 
\cite{givental-mirrorthm-projective}. 
We assume that the $\T$-action on $X$ is Hamiltonian. 
For projective manifolds, 
this assumption is equivalent to 
that the action has at least one fixed point.
Then, the equivariant cohomology is of the form 
$H^*_{\T}(X)\cong H^*(X)\otimes H^*_\T({\rm pt})$. 
This isomorphism is not canonical, but we choose 
homogeneous equivariant lifts $p_0,\dots,p_s$ 
of the basis in section \ref{sect:quantumD-modules}. 
Let $t_0,\dots, t_s$ be linear coordinates dual to 
$p_0,\dots,p_s$ 
and put $H^*_{\T}({\rm pt})=\C[\lambda_1,\dots,\lambda_l]$. 
Set $q_a=\exp(t_a)$ for $1\le a\le r$. 
The equivariant quantum cohomology is of the form
\[
QH^*_{\T}(X)=H^*_{\T}(X)\hatotimes 
\C[\![t_0,q_1,\dots,q_r,t_{r+1},\dots,t_s]\!]
=H^*(X)\otimes\C[\lambda][\![x]\!].
\]
Here, $\C[\lambda][\![x]\!]$ is shorthand for 
$\C[\lambda_1,\dots,\lambda_l][\![t_0,q_1,\dots,q_r,t_{r+1},\dots,t_s]\!]$ 
and $\hatotimes$ means the tensor product 
completed in the $x$-adic topology.

Our $(QH^*_{\T}(X),*)$ defines a 
{\it quasi-conformal Frobenius structure} \cite{givental-elliptic} 
on the base space 
$\mM=\Spec(\C(\lambda)[\![x]\!][q_1^{-1},\dots,q_r^{-1}])$. 
Let $\mT$ be the relative tangent sheaf of 
$\mM\rightarrow \Spec\C(\lambda)$. 
Then, each element of $QH^*_\T(X)$ gives a section of $\mT$ 
by the correspondence $p_i\mapsto\partial_i:=\partial/\partial t_i$. 
The quasi-conformal Frobenius structure of $QH_{\T}^*(X)$
consists of the following data: 

\begin{enumerate}
\item Flat $\C(\lambda)$-bilinear 
symmetric pairing $\pairT{\cdot}{\cdot}$ on $\mT$ 
defined by $\pairT{\partial_i}{ \partial_j}=\int_X^\T p_i\cup p_j$.  

\item $\C(\lambda)$-bilinear symmetric product 
$*\colon \mT\otimes\mT\rightarrow \mT$ 
satisfying $\pairT{\partial_i*\partial_j}{\partial_k}
=\pairT{\partial_i}{\partial_j*\partial_k}$. 

\item Flat unit section $1=p_0=\partial_0\in \mT$. 

\item Euler operator $\tEu$ defined by 
\[
\tEu=2 t_0\parfrac{}{t_0}+2 c_0(\lambda)\parfrac{}{t_0}
+\sum_{a=1}^r 2 c_a \parfrac{}{t_a}
+\sum_{j=r+1}^{s} (2-\deg p_j) t_j\parfrac{}{t_j}+
\sum_{j=1}^l 2 \lambda_j \parfrac{}{\lambda_j},
\]
where 
$c_0(\lambda)+\sum_{a=1}^r c_a p_a=c_1^\T(TX)$. 
\end{enumerate}
The name {\it quasi}-conformal 
comes from that $\tEu$ is {\it not} a section of $\mT$ 
(it contains the derivation $\lambda_j\partial/\partial\lambda_j$). 
As in non-equivariant case, 
the dual Givental connection 
$\nabla_j^\hbar=\hbar(\partial/\partial t_j)+p_j*$ on 
$\mT$ is flat for any value of $\hbar$. 
This defines the $\T$-equivariant quantum $D$-module: 
\[
QDM_\T^*(X)=(H^*(X)\otimes\C[\lambda,\hbar][\![x]\!],\nabla^\hbar). 
\]
The Euler operator satisfies the following: 
\begin{align*}
[\tEu, V_1*V_2]&=V_1*[\tEu,V_2]+[\tEu,V_1]*V_2+2V_1*V_2, \\
\tEu\pairT{V_1}{V_2}&=\pairT{[\tEu,V_1]}{V_2}
+\pairT{V_1}{[\tEu,V_2]}+(4-2N)\pairT{V_1}{V_2}, 
\end{align*}
for $V_1,V_2\in \mT$ and $N=\dim_{\C} X$.

We assume some familiarity with  
the semisimplicity and canonical coordinates 
for conformal Frobenius manifolds, 
see e.g. \cite{dubrovin,manin}. 
Here, we review the construction 
of canonical coordinates in equivariant quantum cohomology. 
Assume that the $\T$ action on $X$ has only isolated fixed points. 
Then, by the localization theorem of equivariant cohomology, 
we have an isomorphism of rings
\[
H^*_{\T}(X)\otimes_{\C[\lambda]}\C(\lambda)
\cong \bigoplus_{\sigma\in X^{\T}} \C(\lambda)
\psi_\sigma, \quad 
\psi_{\sigma}=\frac{i_{\sigma *}([\sigma])}{\be^{\T}(T_\sigma X)}.
\]
Here, $i_{\sigma}\colon \{\sigma\} \rightarrow X$ is the inclusion and 
$\be^\T$ is the $\T$-equivariant Euler class. 
Because $i_\sigma^*(\psi_{\tau})=\delta_{\sigma \tau}$ for 
$\sigma, \tau\in X^\T$, 
we have $\psi_\sigma\cup\psi_\tau=\delta_{\sigma\tau}\psi_\sigma$. 
The localized equivariant quantum cohomology is also semisimple. 
The idempotent $\psi_\sigma$ can be deformed to the 
idempotent $\psi_\sigma(x)$ of $QH^*_\T(X)\hatotimes_{\C[\lambda]}\C(\lambda)$
such that $\lim_{q\to 0}\psi_{\sigma}(x)=\psi_\sigma$ and 
$\psi_\sigma(x)*\psi_\tau(x)=\delta_{\sigma\tau}\psi_\sigma(x)$. 
The projection to the $\psi_\sigma(x)$-component 
$QH^*_{\T}(X)\hatotimes_{\C[\lambda]}\C(\lambda)
\rightarrow \C(\lambda)[\![x]\!]$ corresponds to 
a closed 1-form on $\mM$. 
A primitive $\cu_\sigma(t,\lambda)$ 
of this closed 1-form is called a {\it canonical coordinate}. 
Then, we can write  
$\psi_\sigma(x)=\partial/\partial \cu_\sigma$. 
Note that $d\cu_\sigma$ is characterized by 
\begin{equation}
\label{eq:dofcanonical}
V*\psi_\sigma(x)=d\cu_\sigma(V) \psi_\sigma(x), \quad V\in \mT.
\end{equation}
The canonical coordinates $\{\cu_\sigma\}_{\sigma\in X^\T}$ 
are determined up to functions in $\lambda_i$.

Put $\alpha(\sigma)=i_\sigma^*(\alpha)$ for $\alpha\in H^*_\T(X)$. 
Let $T_\sigma^*X\cong \chi_1(\sigma)\oplus\cdots \oplus\chi_N(\sigma)$
be the weight decomposition, where 
$\chi_j(\sigma)\in\C[\lambda]$. 
Although $\mM$ does not contain the divisor $q_a=0$,  
the limit $\lim_{q\to 0}d\cu_\sigma$ exists 
as an element of $(H^*(X))^*\otimes \C(\lambda)$ 
because the dual basis $\psi_\sigma(x)$ is in 
$H^*(X)\otimes\C(\lambda)[\![x]\!]$.  
Since $\lim_{q\to 0} d\cu_\sigma$ is the projection to 
$\psi_\sigma$, we have  
$d\cu_\sigma = \sum_{i=0}^s (p_i(\sigma)+O(q))dt_i$, 
where $O(q)\in \C(\lambda)[\![x]\!]$.  
We normalize canonical coordinates by the following 
classical limit condition: 
\begin{equation}
\label{eq:normalizationofcanonical}
\cu_\sigma(t,\lambda)=\sum_{j=1}^N (-\chi_j(\sigma) 
+\chi_j(\sigma)\log (-\chi_j(\sigma))) 
+\sum_{i=0}^s p_i(\sigma) t_i +O(q), 
\end{equation}
where $O(q)\in \C(\lambda)[\![x]\!]$. 
This together with (\ref{eq:dofcanonical})
determines $\cu_\sigma(t,\lambda)$ uniquely. 
Using $c_1^\T(T_\sigma X)
=c_0(\lambda)+\sum_{a=1}^r c_ap_a(\sigma)=-\sum_{j=1}^N \chi_j(\sigma)$, 
we can check that this $\cu_\sigma$ automatically satisfies the homogeneity  
$\tEu \cu_\sigma = 2\cu_\sigma$.

Assume that $QH^*_\T(X)$ has convergent structure constants
around $x=\lambda=0$.   
Then $\cu_\sigma(t,\lambda)$ becomes a 
multi-valued analytic function. 
We conjecture the following which resembles the $R$-conjecture. 

\begin{conjecture}
\label{conj:u-conjecture}
When the non-equivariant quantum cohomology is 
generically semisimple, 
canonical coordinates $\cu_\sigma(t,\lambda)$ normalized by 
(\ref{eq:normalizationofcanonical}) 
is regular at $\lambda=0$ for a semisimple point $t$ 
of non-equivariant theory.   
\end{conjecture}

Let $\cu'_\sigma(t,\lambda)$ be any canonical coordinate 
which is homogeneous $\tEu\cu'_\sigma = 2 \cu'_\sigma$ and 
is regular at $\lambda=0$ for semisimple $t$. 
If this conjecture is true, such a canonical coordinate 
is different from the above normalization 
only by a linear form $\sum_{i=1}^l a_i\lambda_i$. 
Note that the non-equivariant limit $\lambda\to 0$ of 
$\cu_\sigma(t,\lambda)$ gives the canonical coordinate 
$\ncu(t)$ of non-equivariant theory satisfying the homogeneity 
$E\ncu(t)=2\ncu(t)$. 
(In non-equivariant theory, homogeneous canonical coordinates are  
unique up to order.) 
Later, we will see that this conjecture (together with 
$R$-conjecture stated below) holds for toric variety
using the equivariant mirror.

Consider the following differential equation for a section 
$s(t,\lambda,\hbar) \in \mT$.   
\begin{equation}
\label{eq:qde}
\hbar\parfrac{}{t_i}s(t,\lambda,\hbar)=p_i*s(t,\lambda,\hbar). 
\end{equation}
An {\it asymptotic solution} 
\cite{dubrovin,givental-elliptic,givental-semisimple,givental-quadratic}
$s_\tau(t,\lambda,\hbar)$ for this differential equation 
is the solution of the form 
\[
s_{\tau}(t,\lambda,\hbar)=\sum_{\sigma} 
 R_{\sigma\tau}(t,\lambda,\hbar) e^{\cu_\tau(t,\lambda)/\hbar}
 \hess_\sigma \parfrac{}{\cu_\sigma}, \quad 
\hess_\sigma=\pairT{\parfrac{}{\cu_\sigma}}
 {\parfrac{}{\cu_\sigma}}^{-1/2}, 
\]
where $R_{\sigma\tau}(t,\lambda,\hbar)=\delta_{\sigma\tau}
+R^{(1)}_{\sigma\tau}(t,\lambda)\hbar
+R^{(2)}_{\sigma\tau}(t,\lambda)\hbar^2+\cdots$
is an asymptotic series in $\hbar$. 
The matrix $R_{\sigma\tau}$ 
is uniquely determined by the classical limit condition: 
\begin{equation}
\label{eq:normalizationofR}
\lim_{q\to 0} R_{\sigma\tau}(t,\lambda,\hbar)=
\delta_{\sigma\tau} \exp(b_\sigma(\hbar,\lambda)), \quad  
b_\sigma(\hbar,\lambda)=\sum_{k=1}^\infty 
N_{2k-1}(\sigma)
\frac{B_{2k}}{2k}\frac{\hbar^{2k-1}}{2k-1},
\end{equation}
where 
$N_{2k-1}(\sigma)=\sum_{j=1}^N 1/\chi_{j}(\sigma)^{2k-1}$. 
Since $\lim_{q\to 0}\hess_\sigma=\be^\T(T_\sigma X)^{1/2}$, 
we can see that $R_{\sigma\tau}^{(n)}(t,\lambda)$ belongs to 
$\C(\lambda,\be^\T(T_\sigma X)^{1/2})[\![x]\!]$. 
Under this normalization, $R_{\sigma\tau}$ automatically satisfies  
the following homogeneity and the unitarity: 
\begin{gather*}
(2\hbar\parfrac{}{\hbar}+\tEu) R_{\sigma\tau}(t,\lambda,\hbar) = 0, \quad  
\sum_{\nu\in X^\T} R_{\nu\sigma}(t,\lambda,-\hbar)R_{\nu\tau}(t,\lambda,\hbar)
=\delta_{\sigma\tau}.
\end{gather*}

\begin{conjecture}[$R$-conjecture {\cite{givental-quadratic}}]
\label{conj:R-conjecture}
When the non-equivariant quantum cohomology is generically semisimple, 
the asymptotic solution $R_{\sigma\tau}^{(n)}(t,\lambda)$ normalized by 
(\ref{eq:normalizationofR}) is regular at $\lambda=0$ 
for a semisimple point $t$ of non-equivariant theory. 
\end{conjecture}

An asymptotic solution  $R_{\sigma\tau}(t,\lambda,\hbar)$ 
satisfying the homogeneity and regularity at $\lambda=0$
for a semisimple $t$ is unique if it exists. 
Thus, if $R$-conjecture is true, such a solution exists and satisfies 
the above classical limit condition (\ref{eq:normalizationofR}). 
In this case, the non-equivariant limit 
$\lim_{\lambda\to 0}R_{\sigma\tau}$ gives the 
homogeneous asymptotic solution of non-equivariant theory.  

By the theory of Givental \cite{givental-semisimple,givental-quadratic}, 
the $R$-conjecture implies the Virasoro constraints for 
the non-equivariant Gromov-Witten theory of $X$. 

\subsection{Equivariant mirror} 
We review an equivariant version of the mirror of toric variety 
\cite{givental-mirrorthm-toric}. 
We use the same notation as in section \ref{sect:mirrorsandtheMellinsystem}. 
The torus $\T:=(S^1)^{r+N}$ acts on a toric variety 
$X=\C^{r+N}_{\mathcal{B}}/(\C^*)^r$ by 
$(z_1,\dots,z_{r+N})\mapsto (t_1z_1,\dots,t_{r+N}z_{r+N})$. 
Here, $\td_1,\dots,\td_r\in H^2_\T(X)$ denotes the $\T$-equivariant 
Poincar\'{e} duals of torus invariant prime divisors 
$D_1,\dots,D_{r+N}$.  
In other words, $w_i$ is the $\T$-equivariant first Chern class of 
the line bundle 
\[
\C^{r+N}_{\mathcal{B}}\times\C/
(z_1,\dots,z_{r+N},v)\sim 
(t^{\td_1}z_1,\dots,t^{\td_{r+N}}z_{r+N}, t^{\td_i}v), 
(t_1,\dots,t_r)\in (\C^*)^r
\]
with $\T$-action 
$(z_1,\dots,z_{r+N},v)\mapsto 
(t_1z_1,\dots,t_{r+N}z_{r+N},t_iv)$. 
Let $p_a$ denote an equivariant lift of the nef integral basis 
in section \ref{sect:mirrorsandtheMellinsystem}. 
More precisely, we define $p_a$ as  
the $\T$-equivariant first Chern class of the 
line bundle 
\[
\C^{r+N}_{\mathcal{B}}\times\C/
(z_1,\dots,z_{r+N},v)\sim 
(t^{\td_1}z_1,\dots,t^{\td_{r+N}}z_{r+N}, t_a v), 
(t_1,\dots,t_r)\in (\C^*)^r. 
\]
with $\T$-action 
$(z_1,\dots,z_{r+N},v)\mapsto (t_1z_1,\dots,t_{r+N}z_{r+N},v)$. 
Then, $\td_i$ is written as 
\[
\td_i=\sum_{a=1}^r m_{ia} p_a - \lambda_j.
\]
Define an equivariant phase function $F^\T(\mc;\lambda)$ on the mirror family 
$\pi\colon Y\rightarrow (\C^*)^r$ by
$F^\T(\mc;\lambda)=\sum_{i=1}^{r+N}(\mc_i+\lambda_i\log \mc_i)$.
Set $F^\T_q=F^\T|_{Y_q}$. 
Then consider the equivariant oscillatory integral: 
\[
\mI_\Gamma(q,\lambda,\hbar)
=\int_{\Gamma_q\subset Y_q} e^{F^\T_q/\hbar} \omega_q. 
\]
We can easily obtain the equivariant version of 
Proposition \ref{prop:noneqmirrordiffeq}. 
\begin{proposition}
\label{prop:equivariantMellindiffeq}
The equivariant oscillatory integral $\mI_{\Gamma}$ satisfies 
$\mPT_\bd\mI_{\Gamma}(q,\lambda,\hbar)=0$ for all $\bd\in \Z^r$, where 
\[
\mPT_\bd=q^\bd \prod_{\pair{\td_i}{\bd}<0} \prod_{k=0}^{-\pair{\td_i}{\bd}-1}
\left(\sum_{a=1}^r m_{ia} \hbar\partial_a -\lambda_i -k\hbar\right)
-\prod_{\pair{\td_i}{\bd}>0} 
\prod_{k=0}^{\pair{\td_i}{\bd}-1}
\left(\sum_{a=1}^r m_{ia} \hbar\partial_a -\lambda_i- k\hbar\right). 
\]
\end{proposition} 
We introduce the $D$-module $M_{\rm Mell}^\T$ as in section 
\ref{sect:mirrorsandtheMellinsystem}.   
\[
M_{\rm Mell}^\T = \C\langle q^{\pm},\hbar\partial,\hbar,\lambda\rangle /
I_{\rm Mell}^\T, \quad 
I_{\rm Mell}^\T=\sum_{\bd\in \Z^r} 
\C\langle q^{\pm},\hbar\partial,\hbar,\lambda\rangle \mPT_\bd
\]
\subsection{$\T\times S^1$-equivariant Floer cohomology}
We can construct $\T\times S^1$-equivariant Floer cohomology 
$\FHT$
in an analogous manner to section \ref{sect:S^1-equivariantFloercohomology}.
First, $H_{\T\times S^1}^{\infty/2}$ is defined by the limit of 
the inductive system 
($H^*_{\T\times S^1}(L_{\bd}^\infty)$, $i_{\bd,\bd'}$), 
where $H^*_{\T\times S^1}(L_{\bd}^\infty)=\C[P,\hbar,\lambda]$, 
$i_{\bd,\bd'}(\alpha)=\alpha\prod_{i=1}^N 
\prod_{k=\pair{\td_i}{\bd'}}^{\pair{\td_i}{\bd}-1}W_{i,k}$ and  
$W_{i,k}=\sum_{a=1}^r m_{ia}P_a-\lambda_i-k\hbar$. 
It has an action of $\C\langle P,Q^{\pm},\hbar,\lambda\rangle$. 
Let $\DeltaT$ be the image of $1$ in $H^*_{\T\times S^1}(L_{0}^\infty)$.
Let $FH_0^\T$ be a $\C\langle P,Q,\hbar,\lambda\rangle$-submodule 
of $H^{\infty/2}_{\T\times S^1}$ 
generated by $\DeltaT$. 
Finally, we define $\FHT$ as the $Q$-adic completion 
of $FH_0^\T$. 
We can check that (by the argument in \cite{iritani-EFC}) 
$\FHT$ is a free module over $\C[\hbar,\lambda][\![Q]\!]$ 
and has $\{T_i(P)\DeltaT\}_{i=0}^s$ as a basis, 
where $T_i(y_1,\dots,y_r)$ is the polynomial in Proposition 
\ref{prop:EFCdefinedovermO}.  
By the same argument as Proposition \ref{prop:H^infty/2=Mell}, 
\ref{prop:FH=completedMell}, we obtain 
\begin{gather*}
H^{\infty/2}_{\T\times S^1}\cong M_{\rm Mell}^\T, 
\quad 
FH_0^\T\cong \C\langle q,\hbar\partial,\hbar,\lambda\rangle /I_{\rm poly}^\T, 
\quad 
\FHT\cong \C[\hbar,\lambda][\![q]\!]\langle \hbar\partial\rangle
\Big / \ov{I_{\rm poly}^\T},
\end{gather*}
where 
$I_{\rm poly}^\T=I_{\rm Mell}^\T\cap\C\langle q,
\hbar\partial,\hbar,\lambda\rangle$ and 
$\ov{I^\T_{\rm poly}}$ is the closure in the $q$-adic topology. 
The relationship between the $\T\times S^1$-equivariant Floer cohomology 
and $\T$-equivariant quantum cohomology is given as follows
(c.f. Theorem \ref{thm:generalizedmirrortransformation}). 
\begin{theorem}
\label{thm:equivariantgenmir}
There exists an embedding 
$\emb\colon (\C^r,0) \rightarrow (\C^{s+1},0)$ 
and an isomorphism of $D$-modules
$\Phi_{\emb}\colon \emb^*(QDM_\T^*(X))\cong \FHT$. 
The map $\emb$ is given by the equations 
\[
\hatt_0=F_0(q;\lambda),
\hq_1=q_1\exp(F_1(q;\lambda)),\dots, 
\hq_r=q_r\exp(F_r(q;\lambda)),  
\hatt_{r+1}=F_{r+1}(q;\lambda),\dots, \hatt_s=F_{s}(q;\lambda)
\]
for  $F_i(q;\lambda)\in \C[\lambda][\![q]\!], \ F_i(0;\lambda)=0$ 
and $\Phi_{\emb}|_{q=0}$ is determined by 
the canonical isomorphism 
$H^*_\T(X)\otimes\C[\hbar]\cong \FHT/\sum_{a=1}^r q_a \FHT$. 
\end{theorem} 
This theorem is a generalization of the result of 
\cite{iritani-genmir} to the $\T$-equivariant case.  
The proof is completely similar and based on a $\T$-equivariant version 
of Coates-Givental's quantum Lefschetz theorem 
\cite{coates-givental}. 
The main points in the proof are
(1) any toric variety can be realized as a complete intersection 
of torus invariant divisors in a Fano toric variety $\widetilde{X}$,
(2) 
the $J$-function of $\FHT$ 
is identical with the Coates-Givental modification of 
the $J$-function of small quantum cohomology 
of the ambient $\widetilde{X}$, 
and 
(3) the reconstruction from small to big is unique. 
We can see that Coates-Givental's quantum Lefschetz theorem admits  
a $\T$-equivariant generalization.

\subsection{Asymptotic solution via equivariant mirror}

Using the equivariant mirror, 
we will construct asymptotic solutions 
subject to the classical limit condition (\ref{eq:normalizationofR}).  
Then it follows that the $R$-conjecture \ref{conj:R-conjecture}
holds for toric varieties. 
Unfortunately, 
we could not prove the convergence of $\T$-{\it equivariant} 
quantum cohomology using the method in this paper. 
Instead, 
we use a result in a forthcoming paper \cite{iritani-localization} 
which proves the convergence of equivariant quantum cohomology  
by the localization method.

We start with the construction of canonical coordinates 
for equivariant mirror. 
Each fixed point $\sigma\in X^\T$ can be written as an intersection 
of $N$ toric divisors $\cap_{i\in I_\sigma} D_i$, 
where $I_\sigma$ is a subset of $\{1,2,\dots,r+N\}$. 
We can use $\{\mc_i\}_{i\in I_\sigma}$ as fiber coordinates 
of the family $\pi\colon Y\rightarrow (\C^*)^r$. 
We write for $j\notin I_\sigma$, 
\[
\mc_j=\prod_{a=1}^r q_a^{l^\sigma_{aj}} \prod_{i\in I_\sigma} 
\mc_i^{x^\sigma_{ji}}.
\]
Here, the matrix $l^\sigma_{aj}$ is the inverse of 
$(m_{ja})_{j\notin I^\sigma, a=1,\dots,r}$
and $x_{ji}^\sigma =-\sum_{a=1}^r m_{ia}l^\sigma_{aj}$. 
Since $\{p_a\}_{a=1}^r$ is a nef basis, 
it follows that the matrix $l^\sigma_{aj}$ satisfies $l^\sigma_{aj}\ge 0$ 
and $\sum_{a=1}^r l^\sigma_{aj}>0$.

We need the following elementary lemma. 
See the Appendix for the proof.  
\begin{lemma}
\label{lem:cps}
Let $U\subset \C^l$ be a neighborhood of $0$ and 
$D\subset U\setminus\{0\}$ be the complement of an analytic subvariety in $U$. 
Let $f(q_1,\dots,q_r,\lambda_1,\dots,\lambda_l)$ be a function 
holomorphic in the neighborhood of $\{q=0,\lambda\in D\}$. 
Assume that $f(\cdot,\lambda)$ can be expanded in the form 
$\sum_{\bd\ge 0} f_\bd (\lambda) q^\bd$ 
for each $\lambda\in D$ and that $f_\bd(\lambda)$ can be analytically 
continued to a holomorphic function on $U$. 
Then $f$ can be extended to a holomorphic function around $q=\lambda=0$. 
\end{lemma} 

Put $\chi_i(\sigma):=-\td_i(\sigma)$. 
Then we have the following relations. 
\begin{gather}
\label{eq:chisigma}
T_\sigma^*X\cong \bigoplus_{i\in I_{\sigma}} \chi_i(\sigma),\quad 
\chi_i(\sigma)=\lambda_i+\sum_{j\notin I_\sigma}
 \lambda_j x^\sigma_{ji} \quad \text{ for }i\in I_\sigma,
\\
\nonumber
\chi_j(\sigma)= 
\lambda_j-\sum_{a=1}^r m_{ja}p_a(\sigma)=0\quad 
\text{ for }j\notin I_\sigma,
\end{gather}

\begin{lemma}
\label{lem:mirrorcanonical}
{\rm (i)} For each $\sigma\in X^\T$, there exists a unique branch of 
critical points $\{\crit_\sigma(q)\}_{q}$ of $F^\T_q$ such that 
in a coordinate system $\{\mc_i\}_{i\in I_\sigma}$, 
\[
\mc_i(\crit_\sigma(q))=-\chi_i(\sigma) +O(q), 
\quad i\in I_\sigma.
\]
Here, $O(q)$ has an expansion of the form 
$\sum_{|\bd|>0}f_\bd(\lambda)q^\bd$ for 
$f_\bd(\lambda)\in\C[\lambda,\chi_i(\sigma)^{-1};i\in I_\sigma]$. 

{\rm (ii)} The critical value $F_q^\T(\crit_q(\sigma))$
equals the pull-back $\emb^*(\cu_\sigma)$ of 
the canonical coordinate satisfying the classical limit condition 
$(\ref{eq:normalizationofcanonical})$, 
where $\emb$ is the map in Theorem \ref{thm:equivariantgenmir}. 

{\rm (iii)} The defining equations of the map $\emb$ 
are convergent power series. 
\end{lemma}
\begin{proof}
(i) follows from a direct calculation. 
For (ii), first note that there exists an isomorphism 
\[
dF^\T\colon \bigcup_{q\in (\C^*)^r} \Crit(F^\T_q) \cong 
\Ch(M_{\rm Mell}^\T)
\]
as in non-equivariant case 
(Lemma \ref{lem:criticalfamilyischaracteristic}). 
Because we have a surjection
\[
M_{\rm Mell}^\T\otimes_{\C[\hbar,\lambda,q^{\pm}]}
\C[\hbar,\lambda][\![q]\!][q^{-1}] \rightarrow 
\FHT[q^{-1}]\cong \emb^*(QDM_{\T}^*(X))[q^{-1}],
\]  
$\emb^*(\Ch(QDM_\T^*(X)))|_{q_a\neq 0}$ is 
embedded in $\bigcup_q \Crit(F_q^\T)$. 
The characteristic variety of $QDM_\T^*(X)$ is 
the union $\bigcup_\sigma \Graph(d\cu_\sigma)$ 
of graphs of the differentials of canonical coordinates.  
Thus, for a branch $\{\crit(q)\}$ contained in  
$\emb^*(\Ch(QDM^*_\T(X)))$,   
the critical value $F^\T_q(\crit(q))$
gives us a pull-back of a canonical coordinate. 
Furthermore, by the relations (\ref{eq:chisigma}), 
it follows that 
$F_q^\T(\crit_\sigma(q))$ satisfies the classical limit condition 
(\ref{eq:normalizationofcanonical}). 
From this we can see that the branch $\{\crit_\sigma(q)\}$ 
actually corresponds to a branch of $\Ch(QDM^*_\T(X))$ 
and $F_q^\T(\crit_\sigma(q))=\emb^*(\cu_\sigma)(q)$. 
In the forthcoming paper \cite{iritani-localization}, we will prove 
that $QH^*_\T(X)$ has convergent structure constants. 
Thus, the natural flat coordinates 
$\hx=\{\hatt_0,\hq_1,\dots,\hq_r,\hatt_{r+1},\dots,\hatt_{s}\}$
of $QH_\T^*(X)$ can be written as analytic functions in 
$\cu_\sigma$ and $\lambda$.  
Because $F_q^\T(\crit_\sigma(q))$ is 
a multi-valued analytic function in $q$ and $\lambda$, 
$\emb^*(\hx_i)$ is a (possibly multi-valued) analytic function in $q$ and $\lambda$. 
It is easy to see that it is holomorphic in the neighborhood of 
$\{q=0,\lambda\text{ generic}\}$ (identically zero at $q=0$). 
Since $\emb^*(\hx_i)$ belongs to $\C[\lambda][\![q]\!]$, 
(iii) follows from Lemma \ref{lem:cps}. 
\end{proof}

\begin{remark}
\rm 
{\rm (i)} 
The above lemma shows that Conjecture \ref{conj:u-conjecture} holds 
for toric varieties. 

{\rm (ii)} 
The branch $\{\crit_{\sigma}(q)\}_q$ in the above lemma 
corresponds to a branch described in Proposition \ref{prop:coveringdegree} {\rm (ii)} 
in the non-equivariant limit.  

{\rm (iii)} Because the map $\emb$ preserves the degree, 
the homogeneity for $\emb^*(\cu_\sigma)=F^\T_q(\crit_\sigma(q))$ 
can be written as
$c_0(\lambda)+
(\sum_{a=1}^r c_a q_a\partial/\partial q_a
+\sum_{j=1}^{r+N}\lambda_j\partial/\partial \lambda_j ) 
F^\T_q(\crit_\sigma(q))=F^\T_q(\crit_\sigma(q))$ 
for $ c_0(\lambda)=-\sum_{j=1}^{r+N}\lambda_j$. 
This can also be shown by a direct calculation. 
\end{remark}

Let $\Gamma(\sigma)$ be the descending Morse cycle of $\Re(F^\T_q/\hbar)$ 
from the critical point $\crit_\sigma(q)$.  
From now, we write $\cu_\sigma$ for $\emb^*(\cu_\sigma)$ 
by abuse of notation.  
Let $\mI_{\sigma}(q,\lambda,\hbar):=\mI_{\Gamma(\sigma)}(q,\lambda,\hbar)$ 
be the oscillatory integral over $\Gamma(\sigma)$. 
By the method of stationary phase, 
$\mI_{\sigma}$
can be expanded in an asymptotic series $\mIas_\sigma$. 
Put $\mc_i=e^{\mcT_i} \mc_i(\crit_\sigma)$ and 
$\mcT_i=\sqrt{\hbar}\mct_i$. 
Note that $\mcT_j=\sum_{i\in I_\sigma} x_{ji}^\sigma \mcT_i$ 
for $j\notin I_\sigma$. 
\begin{align}
\nonumber
\mI_{\sigma}(q,\lambda,\hbar)
&=e^{F^\T_q(\crit_\sigma(q))/\hbar}\int_{\Gamma(\sigma)} 
\exp \frac{F^\T_q(e^\mcT\crit_\sigma(q))-F_q^\T(\crit_\sigma(q))}{\hbar}
\prod_{i\in I_\sigma}d\mcT_i \\ 
\label{eq:toricosci}
&= e^{\cu_\sigma/\hbar}\int_{\Gamma(\sigma)} 
\exp
\frac{\sum_{i=1}^{r+N}(e^{\mcT_i}-1)\mc_i(\crit_\sigma(q))+\lambda_i \mcT_i}
{\hbar} 
\prod_{i\in I_\sigma} d\mcT_i  \\
\nonumber
\sim \hbar^{N/2} & e^{\cu_\sigma/\hbar} 
\int_{\R^N} e^{\sum_{i=1}^{r+N} \mc_i(\crit_\sigma)\mct_i^2/2}
\exp\left (\sum_{i=1}^{r+N}\sum_{n=3}^\infty 
\hbar^{n/2-1}\mc_i(\crit_\sigma)\frac{\mct_i^n}{n!} \right) 
\prod_{i\in I_\sigma} d \mct_i \\
\nonumber
& =:\mIas_\sigma(q,\lambda,\hbar)= (2\pi\hbar)^{N/2} e^{\cu_\sigma/\hbar}
\frac{1}{\sqrt{\pm\Hess_\sigma}}
\sum_{n=0}^\infty \mI_{\sigma n}(q,\lambda)\hbar^n.  
\end{align}
Here,  $\mI_{\sigma 0}=1$, 
$\pm\Hess_\sigma=
(-1)^N\det(\partial^2 F^\T_q/\partial \mcT_i\partial\mcT_j
(\crit_\sigma(q)))_{i,j\in I_\sigma}$ 
and  we assumed $\mc_i(\crit_\sigma)<0$ and $\hbar>0$.   
All half-integer powers of $\hbar$ disappear due to the antisymmetry. 
We can see that 
the analytic function $\mI_{\sigma n}(q,\lambda)$ can be expanded in a $q$-series 
of the form $\sum_\bd f_\bd(\lambda)q^\bd$, where
 $f_\bd(\lambda)\in \C[\lambda,\chi_i(\sigma)^{-1};i\in I_\sigma]$.  

\begin{proposition} 
\label{prop:asymptoticsolution}
The asymptotic series $\mIas_\sigma$ is a formal solution to the $D$-module 
$\FHT$ for each $\sigma \in X^\T$. 
More precisely, for any differential operator 
$f(q,\lambda,\hbar\partial,\hbar)\in \ov{I^{\T}_{\rm poly}}
\subset \C[\hbar,\lambda][\![q]\!]\langle \hbar\partial\rangle$, 
we have $f(q,\lambda,\hbar\partial,\hbar) \mIas_\sigma(q,\lambda,\hbar)=0$ 
as an $\hbar$-series.  
\end{proposition}
\begin{proof}
Take $f\in \ov{I_{\rm poly}^\T}$. 
We put $f\cdot\mIas_\sigma=(2\pi\hbar)^{N/2} e^{\cu_\sigma/\hbar} 
\sum_{k=0}^\infty a_k(q,\lambda)\hbar^k$. 
Note that $a_k(q,\lambda)$ is a formal power series in $q$. 
Let $f_n\in I^\T_{\rm poly}$ be a sequence converging to 
$f\in \ov{I_{\rm poly}^\T}$ such that $g_n=f-f_n=O(q^n)$. 
Then we have 
\begin{equation}
\label{eq:a_kiszero}
g_n\cdot\mIas_\sigma = 
f\cdot \mIas_\sigma=(2\pi\hbar)^{N/2} e^{\cu_\sigma/\hbar} 
\sum_{k=0}^\infty a_k(q,\lambda)\hbar^k
\end{equation}
by Proposition \ref{prop:equivariantMellindiffeq}.
Because we have the following expansions 
\[
(\pm\Hess_\sigma)^{-1/2}=(\prod_{i\in I_\sigma}\chi_i(\sigma))^{-1/2}(1+\sum_{|\bd|>0}
b_\bd(\lambda)q^\bd), \quad
q_a\parfrac{\cu_\sigma}{q_a}=\sum_{\bd\ge 0}c_\bd(\lambda)q^\bd 
\]  
for some $b_\bd(\lambda), c_\bd(\lambda)\in 
\C[\lambda,\chi_i(\sigma)^{-1};i\in I_\sigma]$, 
we can see from (\ref{eq:a_kiszero}) that $a_k(q,\lambda)=O(q^n)$.  
Because this holds for all $n$, $a_k$ must be zero. 
\end{proof} 

In non-equivariant case, each coefficient of 
the asymptotic expansion of $\mI_{\Gamma}$ in $\hbar$ 
is still analytic function in $q$ 
but {\it not} necessarily expanded in $q$-series. 
However, $f(q,\hbar\partial,\hbar) \cdot \mIas_{\Gamma}(q,\hbar)$ 
does make sense for 
$f\in \mO_{\rm small}^\hbar\langle \hbar\partial\rangle$
as an $\hbar$-series. 

\begin{corollary}
\label{cor:asymptoticsolution}
Let $\Gamma$ be the descending Morse cycle of $\Re(F_q)$ from 
a critical point in a branch described in Proposition \ref{prop:coveringdegree} (ii). 
The asymptotic expansion of the non-equivariant mirror oscillatory integral
$\mI_{\Gamma}(q,\hbar)$ in $\hbar$ gives a formal solution to $\tFH$. 
More precisely, any differential operator in $\ov{I_{\rm poly}}\subset 
\mO_{\rm small}^\hbar\langle \hbar\partial\rangle$ annihilates 
the asymptotic expansion of $\mI_\Gamma$ as an $\hbar$-series.  
\end{corollary}

For $v\in QH^*_\T(X)$, there exists a differential operator 
$\mD[v](q,\lambda,\hbar\partial,\hbar)\in 
\C[\hbar,\lambda][\![q]\!]\langle \hbar\partial\rangle$ such that 
$\Phi_{\emb} (v)=\mD[v](Q,\lambda,P,\hbar)\DeltaT$. 
We can write 
$v=\mD[v](q,\lambda, \nabla^\hbar,\hbar)\cdot \Phi_{\emb}^{-1}\DeltaT$, 
where $\nabla_a^\hbar=\nabla^\hbar_{\emb_*(\partial_a)}$.  
Define a formal section $s_\sigma(q,\lambda,\hbar)$ of $\emb^*\mT$ by the formula 
\[
\pairT{v}{(-2\pi\hbar)^{N/2}s_\sigma}=
\mD[v](q,\lambda,\hbar\partial,\hbar) \cdot\mIas_\sigma. 
\]
By Proposition \ref{prop:asymptoticsolution}, the right hand side 
depends only on $v$ and not on a choice of $\mD[v]$. 
Using $\pairT{\partial/\partial \cu_\sigma}
{\partial/\partial \cu_\tau}=\delta_{\sigma\tau}\hess_\sigma^{-2}$, 
we can also write 
\[
s_\sigma(q,\lambda,\hbar)=\frac{1}{(-2\pi\hbar)^{N/2}}\sum_{\tau} \hess_\tau 
(\mD[\psi_\tau(q)]\mIas_\sigma)\hess_\tau \parfrac{}{\cu_\tau}, 
\]
where $\psi_\tau(q)=\psi_\tau(\emb(q))$. 

\begin{theorem}
\label{thm:R-conjecture}
{\rm (i)} The section $s_\sigma$ is 
the pull-back of an asymptotic solution in section 
\ref{sect:equivariantquantumcohomology}
by the map $\emb$.   

{\rm (ii)} This satisfies the classical limit condition (\ref{eq:normalizationofR}). 

{\rm (iii)} The $R$-conjecture \ref{conj:R-conjecture} holds 
for smooth projective toric variety. 
\end{theorem} 
\begin{proof}
(i): For $v\in H^*_\T(X)$ and $\tilde{s}_\sigma=(-2\pi\hbar)^{N/2}s_\sigma$, 
we have 
\begin{align*}
\pairT{v}{\hbar\parti{a} \tilde{s}_\sigma} 
&=\hbar\partial_a\mD[v] \mIas_\sigma 
=\mD[\nabla^\hbar_{\emb_*(\partial_a )} v] \mIas_\sigma \\
&= \pairT{\emb_*(\partial_a ) * v}{\tilde{s}_\sigma} 
  =\pairT{v}{\emb_*(\partial_a ) * \tilde{s}_\sigma}. 
\end{align*}
Therefore, $s_\sigma$ is a solution to the differential equation (\ref{eq:qde}).
We calculate 
\[
\mD[\psi_\tau(q)] \mIas_\sigma = 
(2\pi\hbar)^{N/2} e^{\cu_\sigma/\hbar} 
\left(
\frac{\mD[\psi_\tau(q)](q,\lambda,\partial\cu_\sigma,\hbar=0)} 
{\sqrt{\pm\Hess_\sigma}} +O(\hbar)   \right).  
\]
Let $V(q,\lambda,\hbar\partial,\hbar)$ be a differential operator 
such that 
$V(q,\lambda,\nabla^\hbar,\hbar)\cdot 1=\Phi_{\emb}^{-1}\DeltaT$. 
Then we have 
$\psi_\tau(q)=\mD[\psi_\tau(q)](q,\lambda,\nabla^\hbar,\hbar)
V(q,\lambda,\nabla^\hbar,\hbar) \cdot 1
=\mD[\psi_\tau(q)](q,\lambda,\partial*,0)V(q,\lambda,\partial*,0)1$, 
where $\partial_a*$ is the quantum product by $\emb_*(\partial_a)$. 
By using (\ref{eq:dofcanonical}),  we have 
\begin{align*}
\delta_{\sigma\tau}\psi_\sigma(q)=\psi_\tau(q)*\psi_\sigma(q) &=
\mD[\psi_\tau(q)](q,\lambda,\partial*,0) V(q,\lambda,\partial*,0) 
\psi_\sigma(q) \\
&=\mD[\psi_\tau(q)](q,\lambda,\partial\cu_\sigma,0)
V(q,\lambda,\partial\cu_\sigma,0)\psi_\sigma(q). 
\end{align*}
Therefore, 
\[
\frac{1}{(-2\pi\hbar)^{N/2}} 
\hess_\tau\mD[\psi_\tau(q)]\mIas_\sigma = 
e^{\cu_\sigma/\hbar} 
\left(\frac{\delta_{\sigma\tau}\hess_\sigma}
{V(q,\lambda,\partial\cu_\sigma,0)\sqrt{\Hess_\sigma}} +O(\hbar)\right). 
\]
The equality 
$\hess_\sigma=V(q,\lambda,\partial\cu_\sigma,0)\sqrt{\Hess_\sigma}$
follows from the differential equation 
$\hbar\partial_a s_\sigma=\partial_a*s_\sigma$, 
and the classical limit 
$\lim_{q\to 0} \sqrt{\Hess_\sigma}=
\sqrt{\prod_{i\in I_\sigma} (-\chi_i(\sigma))}
=\lim_{q\to 0} \hess_\sigma$, $V|_{q=0}=1$. 
Therefore, $s_\sigma$ is an asymptotic solution in section 
\ref{sect:equivariantquantumcohomology}. 

(ii): By Lemma \ref{lem:mirrorcanonical} and (\ref{eq:toricosci}), we have 
\begin{align*}
\hess_\tau \mD[\psi_\tau(q) ] \mI_\sigma 
&=\hess_\tau \mD[\psi_\tau(q) ] 
\left(e^{\cu_\sigma/\hbar} 
\int_{\R^N_{\ge 0}} 
e^{-\sum_{i\in I_\sigma}\chi_i(\sigma)(\mcx_i-\log\mcx_i-1)/\hbar+O(q) } 
\prod_{i\in I_\sigma}d\log \mcx_i\right) \\
=e^{\cu_\sigma/\hbar} 
&\int_{\R^N_{\ge 0}} 
\left(
\hess_\tau\mD[\psi_\tau(q)](0,\lambda,\partial \cu_\sigma, \hbar)
\prod_{i\in I_\sigma} e^{-(\mcx_i-1)\chi_i(\sigma)/\hbar} 
\mcx_i^{\chi_i(\sigma)/\hbar-1}
+O(q)\right)  
\prod_{i\in I_\sigma}d \mcx_i \\
=e^{\cu_\sigma/\hbar} 
& \left(\delta_{\sigma\tau}\sqrt{\prod_{i\in I_\sigma} -\chi_i(\sigma)}
\prod_{i\in I_\sigma} 
e^{\chi_i(\sigma)/\hbar} 
\left(\frac{\chi_i(\sigma)}{\hbar}\right )^{-\chi_i(\sigma)/\hbar}
\Gamma\left(\frac{\chi_i(\sigma)}{\hbar}\right)+O(q)\right). 
\end{align*}
We put $\mcx_i=e^{\mcT_i}$ in the first line,   
used $\partial_a\partial_b\cu_\sigma=O(q)$ in the second line 
and $\mD[\psi_\tau(q)](0,\lambda,\partial \cu_\sigma, \hbar)
=\delta_{\tau\sigma}$ in the third line. 
Note that we are assuming $\chi_i(\sigma)>0$. 
The conclusion follows from the asymptotic expansion of 
the Gamma function around infinity (Stirling's formula).  
\[
\log \Gamma(z) \sim \left(z-\frac{1}{2}\right)\log z - z + \frac{1}{2} \log 2\pi 
+\sum_{k=1}^\infty \frac{B_{2k}}{2k(2k-1)z^{2k-1}}, \quad 
\Re(z)>0.
\]

(iii): First note that if the asymptotic solution is regular at $\lambda=0$ 
for one semisimple point $t_0$, then it is regular for any semisimple point. 
(We can use the value at $t_0$ 
as an initial condition for the differential equation (\ref{eq:qde}).) 
Thus, it suffices to check that 
the above solution $s_\sigma$ is regular at $\lambda=0$ 
for a semisimple $q$. 
Let $T_i(y_1,\dots,y_r)$ be the polynomial in Proposition 
\ref{prop:EFCdefinedovermO}.  
Because 
$\{T_i(P)\DeltaT\}_{i=0}^s$ forms a $\C[\hbar,\lambda][\![q]\!]$-basis 
of $\FHT$, 
we can write 
$\Phi_{\emb}(\partial/\partial t_i)=
\sum_{j=0}^s A_{ij}(q,\lambda,\hbar) T_j(P)\DeltaT$ 
for some $A_{ij}\in \C[\hbar,\lambda][\![q]\!]$.  
Let $\hess_\tau\psi_\tau(q)=
\sum_{i=0}^s \Psi_{\tau i}(q,\lambda) (\partial/\partial t_i)$. 
Because $QH^*_\T(X)$ is convergent, 
$\Psi_{\tau i}(q,\lambda)$ is an analytic function 
of $q,\lambda$ which is regular at $\lambda=0$ for semisimple $q$. 
We can take $\mD[\psi_\tau(q)]$ as 
\[
\hess_\tau\mD[\psi_{\tau}(q)]=
\sum_{i,j=0}^s \Psi_{\tau i}(q,\lambda) 
A_{ij}(q,\lambda,\hbar) T_j(\hbar\partial).
\]
Then, the matrix $R$ is given by 
$(-2\pi\hbar)^{N/2} e^{\cu_\sigma/\hbar} R_{\tau\sigma} =
\sum_{i,j=0}^s  \Psi_{\tau i}(q,\lambda) A_{ij}(q,\lambda,\hbar) 
T_j(\hbar\partial)\mIas_\sigma$. 
Because we can write $T_j(\hbar\partial)\mIas_\sigma = 
(-2\pi\hbar)^{N/2} \mR_{j\sigma} e^{\cu_\sigma/\hbar}$ 
for some $\mR_{j\sigma}=\sum_{n\ge 0}\mR_{j\sigma}^{(n)}(q,\lambda)\hbar^n$,  
we have 
\[
R_{\tau\sigma}(q,\lambda,\hbar)
=\sum_{i,j=0}^s 
\Psi_{\tau i}(q,\lambda) 
A_{ij}(q,\lambda,\hbar) \mR_{j\sigma}(q,\lambda,\hbar).
\]
On the right hand side, 
the regularity at $\lambda=0$ is clear for 
$\Psi_{\tau i}$ and  $\mR_{j\sigma}$. 
Thus, it suffices to show that 
$A_{ij}^{(n)}(q,\lambda)$ is a holomorphic function around $q=\lambda=0$  
where $A_{ij}=\sum_{n\ge 0} A_{ij}^{(n)}(q,\lambda)\hbar^n$. 
Since 
$(A_{ij})=(\Psi_{\tau i})^{-1} (R_{\tau\sigma})(\mR_{j\sigma})^{-1}$ 
and $R_{\tau\sigma}^{(n)}(q,\lambda)$ is an analytic function 
(because $QH^*_\T(X)$ is convergent), 
we can see that $A_{ij}^{(n)}(q,\lambda)$ is convergent around $q=\lambda=0$ 
from Lemma \ref{lem:cps}. 
\end{proof}

\section{Appendix}

\subsection{ Proof of Lemma \ref{lem:epsilongoestozero}} 
First we study $\varepsilon_1(M)$. 
We have  
\begin{align*}
\varepsilon_1(d,M)&:=\sum_{i=1}^{d-1}\frac{(d+1)^M}{(d-i+1)^M(i+1)^M}  \\
&\le \int_{0}^d\frac{(d+1)^M dx}{(d-x+1)^M(x+1)^M}  
=\epsilon^{M-1}(1-\epsilon)^{M}\int_{\epsilon}^{1-\epsilon}
\frac{dy}{y^M(1-y)^M},
\end{align*}
where we put $y=(x+1)/(d+2)$ and $\epsilon=1/(d+2)$. 
Set $y=\sin^2(\theta/2)$ and $\epsilon=\sin^2(\alpha/2)$. 
Then we have 
\[
\varepsilon_1(d,M)\le 4\cos^2\left(\frac{\alpha}{2}\right)
(\sin\alpha)^{2M-2}\int_{\alpha}^{\pi/2}\frac{d\theta}{(\sin\theta)^{2M-1}}. 
\]
Because we can assume $d\ge 2$, we have 
$0<\alpha\le \pi/3$. 
Using the formula (see \cite{calculus}, 2.515)
\[
\int \frac{d\theta}{(\sin\theta)^{2M-1}}
=\frac{(2M-3)!!}{(2M-2)!!}
\left\{-\sum_{r=0}^{M-2}\frac{(2M-4-2r)!!}{(2M-3-2r)!!}
\frac{\cos\theta}{(\sin\theta)^{2M-2-2r}}+
\log\left|\tan \frac{\theta}{2}\right|\right\},
\]
we have 
\begin{align*}
\varepsilon_1(d,M) & \le 4\cos^2\left(\frac{\alpha}{2}\right)
\frac{(2M-3)!!}{(2M-2)!!}
\Biggl\{\sum_{r=0}^{M-2}\frac{(2M-4-2r)!!}{(2M-3-2r)!!}
\cos\alpha(\sin\alpha)^{2r} \\
&\qquad\qquad\qquad  -(\sin\alpha)^{2M-2}
\log\left|\tan \frac{\alpha}{2}\right|\Biggr\} \\
&\le 4 \frac{(2M-3)!!}{(2M-2)!!}
\left\{\frac{1}{1-\sin^2\alpha}+(\sin\alpha)^{2M-3}\frac{2}{e}\right\} \\
&\le 4\left(1-\frac{1}{2M-2}\right)
\left(1-\frac{1}{2M-4}\right)\cdots
\left(1-\frac{1}{2}\right)\left(4+\frac{2}{e}\right). 
\end{align*}
Therefore, we have $\varepsilon_1(M)\to 0$ as $M\to \infty$.

Next we study $\varepsilon_2(B_2,M)$. 
Set $C=B_2/B_1$. We have for $C>1$, 
\begin{align*}
\sum_{i=1}^d \left(\frac{B_1}{B_2}\right)^i\frac{(d+1)^M}{(d+1-i)^M} 
&\le \int_0^d \frac{(d+1)^M dx}{C^x(d+1-x)^M}+\frac{(d+1)^M}{C^d} \\
&\le \int_0^{1-\frac{1}{d+1}} \frac{(d+1)dy}{C^{(d+1)y}(1-y)^M}
+\left(\frac{2M}{e\log C}\right)^M \\
&\le \int_0^{1/2} \frac{(d+1)2^M}{C^{(d+1)y}}dy 
 + \int_{1/2}^{1-\frac{1}{d+1}}\frac{(d+1)dy}{C^{(d+1)/2}(1-y)^M}
 +\left(\frac{2M}{e\log C}\right)^M \\
&\le \frac{2^M}{\log C} + \frac{(d+1)^M}{(M-1)C^{(d+1)/2}}
 +\left(\frac{2M}{e\log C}\right)^M \\
&\le \frac{2^M}{\log C} + \frac{1}{M-1}\left(\frac{2M}{e\log C}\right)^M
 +\left(\frac{2M}{e\log C}\right)^M.  
\end{align*}
Hence we have $\varepsilon(B_2,N)\to 0$ as $B_2\to \infty$. 

\subsection{Proof of Lemma \ref{lem:cps}} 
We prove the lemma by induction on $l$. 
We write $\lambda=(\lambda',\lambda_l)\in \C^{l-1}\times\C$. 
Set $D'=\{\lambda'\in \C^{l-1} \;|\; 
\{\lambda'\}\times\C\cap D \text{ is non-empty}\}$. 
We will show that $f$ is holomorphic in the neighborhood of 
$\{q=0,\lambda_l=0, \lambda'\in D'\}$. 
We choose $\epsilon(\lambda')>0$ so that 
$\{\lambda'\}\times \{|\lambda_l|=\epsilon(\lambda')\}\subset D$. 
Then there exists a positive $C(\lambda')>0$ such that 
$|f_\bd(\lambda',\lambda_l)|\le C(\lambda')^{|\bd|+1}$ 
on $|\lambda_l|=\epsilon(\lambda')$. 
For $|\zeta|<\epsilon(\lambda')/2$, because $f_\bd(\lambda)$ is 
holomorphic on $U$,  
\[
f_\bd(\lambda',\zeta)=\frac{1}{2\pi i} \int_{|z|=\epsilon(\lambda')}
\frac{f_\bd(\lambda',z)}{z-\zeta} dz. 
\] 
From this we can see that $f_\bd(\lambda',\zeta)\le 2 C(\lambda')^{|\bd|+1}$ for  
$\lambda'\in D'$ and $|\zeta|<\epsilon(\lambda')/2$.

\end{document}